\documentclass[a4paper, reqno, 11pt]{amsart}

\usepackage{amscd, amssymb, amsthm, amsmath, enumerate}
\usepackage{graphicx, psfrag, color}

\newtheorem*{MainA}{Main Theorem A}
\newtheorem*{MainB}{Main Theorem B}

\newtheorem{lem}{Lemma}[section]
\newtheorem{thm}[lem]{Theorem}

\newtheorem{cor}[lem]{Corollary}

\newtheorem{defi}[lem]{Definition}

\newcommand{\Stat}{{\mathsf{Stat}}}
\newcommand{\stat}{{\mathsf{stat}}}

\newcommand{\inv}{{\mathsf{inv}}}

\newcommand{\sor}{{\mathsf{sor}}}
\newcommand{\Cyc}{{\mathsf{Cyc}}}
\newcommand{\cyc}{{\mathsf{cyc}}}
\newcommand{\Rmip}{{\mathsf{Rmip}}}
\newcommand{\Rmil}{{\mathsf{Rmil}}}
\newcommand{\rmin}{{\mathsf{rmin}}}

\newcommand{\Lmil}{{\mathsf{Lmil}}}
\newcommand{\lmin}{{\mathsf{lmin}}}
\newcommand{\Lmap}{{\mathsf{Lmap}}}
\newcommand{\Lmal}{{\mathsf{Lmal}}}
\newcommand{\lmax}{{\mathsf{lmax}}}

\newcommand{\Lmic}{{\mathsf{Lmic}}}
\newcommand{\lmic}{{\mathsf{lmic}}}

\newcommand{\Leh}{{\mathsf{Leh}}}
\newcommand{\Max}{{\mathsf{Max}}}
\newcommand{\Min}{{\mathsf{Min}}}

\def\ff{\mathbf{f}}

\def\D{\mathcal{D}}
\def\G{\mathsf{G}}

\numberwithin{equation}{section}

\begin{document}
\title{The sorting index on colored permutations and even-signed permutations}

\author{Sen-Peng Eu}
\address{Department of Mathematics \\
National Taiwan Normal University \\
Taipei, Taiwan 116, ROC}
\email[Sen-Peng Eu]{speu@math.ntnu.edu.tw}

\author{Yuan-Hsun Lo}
\address{Department of Mathematics \\
National Taiwan Normal University \\
Taipei, Taiwan 116, ROC}
\email[Yuan-Hsun Lo]{yhlo0830@gmail.com}

\author{Tsai-Lien Wong}
\address{Department of Applied Mathematics \\
National Sun Yat-sen University \\
Kaohsiung, Taiwan 804, ROC} \email[Tsai-Lien
Wong]{tlwong@math.nsysu.edu.tw}



\thanks{Partially supported by National Science Council, Taiwan under grants NSC 101-2115-M-390-004-MY3 (S.-P. Eu and Y.-H. Lo) and 102-2115-M-110-006-MY2 (T.-L. Wong).}


\maketitle


\begin{abstract}
We define a new statistic $\sor$ on the set of colored permutations
$\G_{r,n}$ and prove that it has the same distribution as the length
function. For the set of restricted colored permutations
corresponding to the arrangements of $n$ non-attacking rooks on a
fixed Ferrers shape we show that the following two sequences of
set-valued statistics are joint equidistributed:
$(\ell,\Rmil^0,\Rmil^1,\ldots,\Rmil^{r-1},$
$\Lmil^0,\Lmil^1,\ldots,\Lmil^{r-1},$
$\Lmal^0,\Lmal^1,\ldots,\Lmal^{r-1},$
$\Lmap^0,\Lmap^1,\ldots,\Lmap^{r-1})$ and
$(\sor,\Cyc^0,\Cyc^{r-1},\ldots,\Cyc^{1}$,
$\Lmic^0,\Lmic^{r-1},\ldots,\Lmic^{1}$,
$\Lmal^0,\Lmal^1,\ldots,\Lmal^{r-1}$,
$\Lmap^0,\Lmap^1,\ldots,\Lmap^{r-1})$. Analogous results are also
obtained for Coxeter group of type $D$. Our work generalizes recent
results of Petersen, Chen-Gong-Guo and Poznanovi\'{c}.
\end{abstract}

\section{Introduction}\label{sec:intro}
\subsection{Mahonian and Stirling statistics}
Let $\mathfrak{S}_n$ be the group of permutations on $n$ letters
$[n]:=\{1,2,\ldots,n\}$. A pair $(\sigma_i,\sigma_j)$ is called an
\emph{inversion} in a permutation $\sigma=\sigma_1\cdots\sigma_n \in
\mathfrak{S}_n$ if $i>j$ and $\sigma_i<\sigma_j$. Denote by
$\inv(\sigma)$ the number of inversions in $\sigma$. The
distribution of $\inv$ over $\mathfrak{S}_n$ was first found by
Rodriguez~\cite{Rodriguez_39} to be
\begin{equation}\label{inversion_distribution}
\sum_{\sigma\in\mathfrak{S}_n}q^{\inv(\sigma)}=\prod_{i=1}^n
[i]_q,
\end{equation}
 where $[i]_q:=1+q+\cdots+q^{i-1}$.

In a Coxeter group, the \emph{length} $\ell(\sigma)$ of a group
element $\sigma$ is the minimal number of generators needed to
express $\sigma$. It is well known~\cite[Chapter 8]{Brenti_05} that
$\mathfrak{S}_n$ is the Coxeter group of type A, where the
generators are the adjacent transpositions and
$\ell(\sigma)=\inv(\sigma)$. A permutation statistic is called
\emph{Mahonian} if it is equidistributed with $\inv$ over
$\mathfrak{S}_n$. Similarly in a Coxeter group a statistic is called
\emph{Mahonian} if it is equidistributed with the length function
$\ell$.

The number of cycles $\cyc$ is another important statistic, whose distribution over $\mathfrak{S}_n$ is  \cite[Proposition 1.3.4]{Stanley_97}
\begin{equation}\label{cycle_distribution}
\sum_{\sigma\in\mathfrak{S}_n}t^{\cyc(\sigma)} =
\prod_{i=1}^n(t+i-1).
\end{equation}
As the coefficients of this polynomial are the (unsigned) Stirling
numbers of the first kind, a permutation statistic over
$\mathfrak{S}_n$ is called \emph{Stirling} if it is equidistributed
with $\cyc$.

The \emph{reflection length} $\ell'(\sigma)$ of $\sigma$ in a Coxeter group is the minimal number of reflections (i.e., elements conjugate to generators) needed to express $\sigma$. 
In type A, the reflections are the transpositions and one has
\begin{equation}\label{eq:cyc_ell'}
\cyc(\sigma)=n-\ell'(\sigma).
\end{equation}

\subsection{Sorting index}
Petersen~\cite{Petersen_11} defined the \emph{sorting index} $\sor$ over $\mathfrak{S}_n$ and proved it is Mahonian. 
On can uniquely decompose $\sigma\in \mathfrak{S}_n$ into a product of transpositions $$\sigma=(i_1\,j_1)(i_2\,j_2)\cdots (i_k\,j_k)$$ with $j_1<j_2<\dots <j_k$ and $i_1<j_1, i_2<j_2, \dots , i_k<j_k$. 
Then the sorting index of $\sigma$ is $$\sor(\sigma):=\sum_{r=1}^k(j_r-i_r).$$
Simply put, $\sor(\sigma)$ counts the number of steps needed to bubble sort a permutation (i.e., the total number of steps needed to successively move $n,n-1, \dots ,1$ back in places). 
For example, for $\sigma=31524=(1\,2)(1\,3)(3\,4)(3\,5)$, the sorting process is
$$31\textbf{5}24\xrightarrow{(3\,5)} 31\textbf{4}25 \xrightarrow{(3\,4)} \textbf{3}1245 \xrightarrow{(1\,3)} \textbf{2}1345 \xrightarrow{(1\,2)} 12345,$$
and $\sor(31524)=(5-3)+(4-3)+(3-1)+(2-1)=6$ as it needs $2,1,2,1$ step(s) respectively to move $5,4,3,2$ back to its place.

By defining $\rmin(\sigma)$ to be the number of right-to-left minima
of $\sigma=\sigma_1\cdots\sigma_n$, Petersen also showed that
$(\inv,\rmin)$ and $(\sor,\cyc)$ have the same joint
distribution~\cite{Petersen_11} and
\begin{equation}~\label{eq:Petersen}
\sum_{\sigma\in\mathfrak{S}_n}q^{\inv(\sigma)}t^{\rmin(\sigma)}
= \sum_{\sigma\in\mathfrak{S}_n}q^{\sor(\sigma)}t^{\cyc(\sigma)}
= \prod_{i=1}^n(t+[i]_q-1).
\end{equation}

The sorting indices $\sor_B$ (for type $B$) and $\sor_D$ (for type
$D$) were also defined by Petersen, and type $B, D$ analogous
identities of ~(\ref{eq:Petersen}) were found by
Petersen~\cite{Petersen_11} and Chen-Gong-Guo~\cite{Chen_13}
respectively. In the case of type $B$, set-valued equidistribution
results are also obtained in~\cite{Chen_13}.

\subsection{Sorting index on a Ferrers shape}
Recently the above results were extended by Poznanovi\'c over the
permutations corresponding to arrangements of $n$ non-attacking
rooks on a fixed Ferrers shape with $n$ rows and $n$
columns~\cite{Poznanovic_14}. The setting can be described as
follows. For a given sequence of integers $\ff=(f_1,f_2,\ldots,f_n)$
with $1\leq f_1\leq f_2\leq \cdots\leq f_n\leq n$, we define the set
of restricted permutations by
$$\mathfrak{S}_{n,\ff} := \{\pi\in\mathfrak{S}_n:\,\pi(i)\leq f_i,
1\leq i\leq n\}.$$

It is clear that $\ff$ defines a Ferrers shape and
$\mathfrak{S}_{n,\ff}$ consists of those permutations corresponding
to non-attacking rook placements. For example, let $\ff=(2,3,3,4)$,
then $\mathfrak{S}_{4,\ff}= \{1234, 1324, 2134, 2314\},$ as
illustrated in Figure 1.
\begin{figure}[h]\label{ferrerperm}
\includegraphics[width=3.5in]{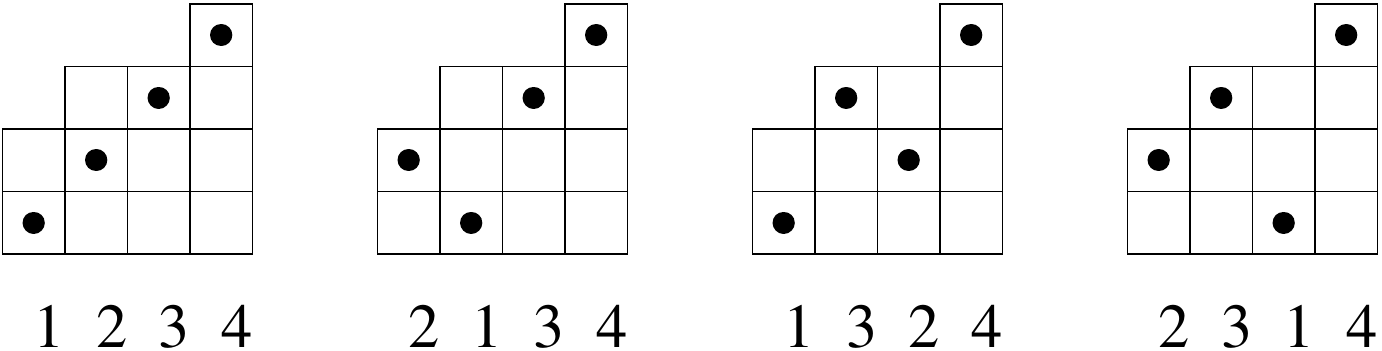}
\caption{The $4$ permutations in $\mathfrak{S}_{4, (2,3,3,4)}$}
\end{figure}
By defining the set-valued statistics
\begin{align*}
\Cyc(\sigma)&:=\{\text{the smallest number in each cycle of the cycle decompositiong}\},\\
\Rmil(\sigma)&:=\{\sigma_i:\,\sigma_i<\sigma_j \text{ for any }j>i\} \mbox{ (Right-to-left minimum letters)},\\
\Lmal(\sigma)&:=\{\sigma_i:\,\sigma_i>\sigma_j \text{ for any }j<i\}\mbox{ (Left-to-right maximum letters)},\\
\Lmap(\sigma)&:=\{i:\,\sigma_i>\sigma_j \text{ for any }j<i\} \mbox{ (Left-to-right maximum places)},
\end{align*}
Poznanovi\'{c}~\cite{Poznanovic_14} proved that
$(\inv,\Rmil,\Lmal,\Lmap)$ and $(\sor,\Cyc,\Lmal,\Lmap)$ have the
same joint distribution over $\mathfrak{S}_{n,\ff}$ by means of
Foata-Han's bijection in \cite{Foata_09}. Analogous results on
Coxeter groups of type $B$ and $D$ were also obtained
in~\cite{Poznanovic_14}, generalizing the works of
Petersen~\cite{Petersen_11} and Chen-Guo-Gong~\cite{Chen_13}.
\medskip

In this paper we extend Ponzanovi\'c's results further in two ways.
In the first part we obtain analogous new results on the colored
permutations $\G_{r,n}$, and in the second part we refine known
results for type $D$.

\subsection{Colored permutations on a Ferrers shape}
The first part of our work is to generalize results on Coxeter group
of type $A(=C_1\wr\mathfrak{S}_n)$ and $B(=C_2\wr\mathfrak{S}_n)$ to
the colored permutations $\G_{r,n}:=C_r\wr\mathfrak{S}_n$ within a
fixed Ferrers shape.

First of all, in Section 2 we will define the sorting index $\sor$
on $\G_{r,n}$ and prove that it is equidistributed with the length
function $\ell$ (defined in the next section).

\begin{thm}[Theorem~\ref{thm:sor_GF}]
For any $r$ and $n$ the statistics $\ell$ and $\sor$ have the same
distribution over $\G_{r,n}$. That is,
\begin{equation}
\sum_{\pi\in\G_{r,n}}q^{\sor(\pi)} = \sum_{\pi\in\G_{r,n}}q^{\ell(\pi)}
=  [n]_q!\cdot \prod_{i=1}^n(1+q^i[r-1]_q).
\end{equation}
\end{thm}

We then consider those colored permutations on a fixed Ferrers shape and seek for analogous results of (\ref{eq:Petersen}) and set-valued equidistributions. 
Denote by $\G_{r,n,\ff}$ the restricted version of $\G_{r,n}$ determined by $\ff$ and define the set-valued statistics $\Cyc^t$, $\Lmic^t$, $\Rmil^t$, $\Lmil^t$, $\Lmal^t$ and $\Lmap^t$ for $t=0,1,\ldots,r-1$ (see Section 2 for detailed definition). 
Our first main theoerm gives two rather interesting long tuples of joint equidistributed set-valued statistics over $\G_{r,n,\ff}$.

\begin{MainA}[Theorem~\ref{thm:equi_restricted}] The two tuples of set-valued statistics
\begin{align*}
\big(&\ell,\Rmil^0,\Rmil^1,\ldots,\Rmil^{r-1},\Lmil^0,\Lmil^1,\ldots,\Lmil^{r-1}, \\
&\hspace{1cm} \Lmal^0,\Lmal^1,\ldots,\Lmal^{r-1},\Lmap^0,\Lmap^1,\ldots,\Lmap^{r-1}\big) \\
\mbox{and}\\
\big(&\sor,\Cyc^0,\Cyc^{r-1},\ldots,\Cyc^{1},\Lmic^0,\Lmic^{r-1},\ldots,\Lmic^{1}, \\
&\hspace{1cm} \Lmal^0,\Lmal^{r-1},\ldots,\Lmal^{1},\Lmap^0,\Lmap^{r-1},\ldots,\Lmap^{1}\big)
\end{align*}
have the same joint distribution over $\G_{r,n,\ff}$.
\end{MainA}
We emphasize that the statistics are set-valued, and the
equidistribution is `twisted', that is, the superindices for each
kind of statistic are $0,1,\dots ,r, r-1$ in the first tuple and
$0,r,r-1,\dots 1$ in the second.

Our second main theorem gives the generating function for counting
colored permutations by the set-valued statistics $(\ell,
\Rmil^0,\Rmil^1,\ldots,\Rmil^{r-1},\Lmil^0,\Lmil^1,\ldots,\Lmil^{r-1})$
or
$(\sor,\Cyc^0,\Cyc^{r-1},\ldots,\Cyc^{1},\Lmic^0,\Lmic^{r-1},\ldots,\Lmic^{1})$.
To understand the statement we need some notation. Given $\ff$ we
define $$H(\ff):=(h_1,h_2,\ldots,h_n),$$ where $h_i$ is the smallest
possible index at which the letter $i$ can appear in a colored
permutation $\sigma\in \G_{r,n,\ff}$. In the above example
$\ff=(2,3,3,4)$ and hence $H(\ff)=(1,1,2,4)$. The superindex of a
statistic is taken as $\bmod \, r$, for example,
$\Rmil^{-2}(\sigma)=\Rmil^{3}(\sigma)$ if $r=5$. The \emph{checking
function} $\xi_{\mathsf{A}}(\omega)$ is defined by
$$
\xi_{\mathsf{A}}(\omega) := \begin{cases}
\omega & \text{if the statement } \mathsf{A} \text{ is true} \\
1 & \text{if the statement } \mathsf{A} \text{ is false}.
\end{cases}
$$

\begin{MainB}[Theorem~\ref{thm:GF_restricted}]
Given $r,n$ and $\ff$. Let
$H(\ff)=(h_1,\ldots,h_n)$. We have
\begin{align*}
&\sum_{\pi\in\G_{r,n,\ff}}q^{\ell(\pi)}\prod_{t=0}^{r-1}\left(\prod_{i\in\Rmil^{-t}(\pi)}x_{t,i}\prod_{i\in\Lmil^{-t}(\pi)}y_{t,i}\right) = \sum_{\pi\in\G_{r,n,\ff}}q^{\sor(\pi)}\prod_{t=0}^{r-1}\left(\prod_{i\in\Cyc^t(\pi)}x_{t,i}\prod_{i\in\Lmic^t(\pi)}y_{t,i}\right) \\
=& ~~\left.\prod_{j=1}^n \right( x_{0,j}+q+\cdots+q^{j-h_j-1}+ \xi_{h_j=1}(y_{0,j})q^{j-h_j} \\
&\qquad + \left.\sum_{t=1}^{r-1} \Big(x_{r-t,j}q^{2j+t-2}+q^{2j+t-3}+\cdots+q^{j+h_j+t-1}+ \xi_{h_j=1}(y_{r-t,j})q^{j+h_j+t-2} \Big) \right).
\end{align*}
\end{MainB}

\subsection{An example}
We illustrate our main theorems by an example. The set
\begin{align*}
\G_{3,2}=\{1^{[0]}2^{[0]}, 1^{[0]}2^{[1]}, 1^{[0]}2^{[2]} , 1^{[1]}2^{[0]}, 1^{[1]}2^{[1]}, 1^{[1]}2^{[2]} , 1^{[2]}2^{[0]} , 1^{[2]}2^{[1]} , 1^{[2]}2^{[2]}, \\
2^{[0]}1^{[0]}, 2^{[0]}1^{[1]}, 2^{[0]}1^{[2]} , 2^{[1]}1^{[0]}, 2^{[1]}1^{[1]}, 2^{[1]}1^{[2]} , 2^{[2]}1^{[0]} , 2^{[2]}1^{[1]} , 2^{[2]}1^{[2]}\},
\end{align*}
consists of $18$ colored permutations of the form
$\sigma_1^{[z_1]}\,\sigma_2^{[z_2]}$, where $\sigma_1\sigma_2$ is a
permutation of $[2]$ and $0\le z_1, z_2 \le 2$ are the colors. The
Main Theorem A says that $$(\ell,\Rmil^0,\Rmil^1,
\Rmil^{2},\Lmil^0,\Lmil^1,\Lmil^{2},
\Lmal^0,\Lmal^1,\Lmal^{2},\Lmap^0,\Lmap^1,\Lmap^{2})$$ and
$$(\sor,\Cyc^0,\Cyc^{2},\Cyc^{1},\Lmic^0,\Lmic^{2},\Lmic^{1},
\Lmal^0,\Lmal^{2},\Lmal^{1},\Lmap^0,\Lmap^{2},\Lmap^{1})$$ have the
same joint distribution, as we can see from Table 1 and Table 2.

Moreover, if we let $\ff=(1,2)$, then $\G_{3,2,\ff}$ consists of the
first $9$ colored permutations and Main Theorem A still holds.

\begin{tiny}
\begin{table}[h]\label{tb1}
\begin{tabular}{c||c|c|c|c|c|c|c|c|c|c|c|c|c}
$\pi$ & $\ell$ & $\Rmil^0$ & $\Rmil^1$ & $\Rmil^2$ & $\Lmil^0$ & $\Lmil^1$ & $\Lmil^2$ & $\Lmal^0$ & $\Lmal^1$ & $\Lmal^2$ & $\Lmap^0$ & $\Lmap^1$ & $\Lmap^2$ \\ \hline \hline
$1^{[0]}2^{[0]}$ & 0 & 1,2 &     &     & 1   &     &     & 1,2 &     &     & 1,2 &     &     \\ \hline
$1^{[0]}2^{[1]}$ & 3 & 1   &   2 &     & 1   &     &     & 1   &   2 &     & 1   &   2 &     \\ \hline
$1^{[0]}2^{[2]}$ & 4 & 1   &     &   2 & 1   &     &     & 1   &     &   2 & 1   &     &   2 \\ \hline
$1^{[1]}2^{[0]}$ & 1 &   2 & 1   &     &     & 1   &     &   2 & 1   &     &   2 & 1   &     \\ \hline
$1^{[1]}2^{[1]}$ & 4 &     & 1,2 &     &     & 1   &     &     & 1,2 &     &     & 1,2 &     \\ \hline
$1^{[1]}2^{[2]}$ & 5 &     & 1   &   2 &     & 1   &     &     & 1   &   2 &     & 1   &   2 \\ \hline
$1^{[2]}2^{[0]}$ & 2 &   2 &     & 1   &     &     & 1   &   2 &     & 1   &   2 &     & 1   \\ \hline
$1^{[2]}2^{[1]}$ & 5 &     &   2 & 1   &     &     & 1   &     &   2 & 1   &     &   2 & 1   \\ \hline
$1^{[2]}2^{[2]}$ & 6 &     &     & 1,2 &     &     & 1   &     &     & 1,2 &     &     & 1,2 \\ \hline \hline
$2^{[0]}1^{[0]}$ & 1 & 1   &     &     & 1,2 &     &     &   2 &     &     & 1   &     &     \\ \hline
$2^{[0]}1^{[1]}$ & 2 &     & 1   &     &   2 & 1   &     &   2 &     &     & 1   &     &     \\ \hline
$2^{[0]}1^{[2]}$ & 3 &     &     & 1   &   2 &     & 1   &   2 &     &     & 1   &     &     \\ \hline
$2^{[1]}1^{[0]}$ & 2 & 1   &     &     & 1   &   2 &     &     &   2 &     &     & 1   &     \\ \hline
$2^{[1]}1^{[1]}$ & 3 &     & 1   &     &     & 1,2 &     &     &   2 &     &     & 1   &     \\ \hline
$2^{[1]}1^{[2]}$ & 4 &     &     & 1   &     &   2 & 1   &     &   2 &     &     & 1   &     \\ \hline
$2^{[2]}1^{[0]}$ & 3 & 1   &     &     & 1   &     &   2 &     &     &   2 &     &     & 1   \\ \hline
$2^{[2]}1^{[1]}$ & 4 &     & 1   &     &     & 1   &   2 &     &     &   2 &     &     & 1   \\ \hline
$2^{[2]}1^{[2]}$ & 5 &     &     & 1   &     &     & 1,2 &     &     &   2 &     &     & 1   \\ \hline
\end{tabular}
\caption{$(\ell,\Rmil^t,\Lmil^t,\Lmal^t,\Lmap^t)$ for $\G_{3,2}$}
\end{table}
\end{tiny}

\begin{tiny}
\begin{table}[h]\label{tb2}
\begin{tabular}{c||c|c|c|c|c|c|c|c|c|c|c|c|c}
$\pi$ & $\sor$ & $\Cyc^0$ & $\Cyc^2$ & $\Cyc^1$ & $\Lmic^0$ & $\Lmic^2$ & $\Lmic^1$ & $\Lmal^0$ & $\Lmal^2$ & $\Lmal^1$ & $\Lmap^0$ & $\Lmap^2$ & $\Lmap^1$ \\ \hline \hline
$1^{[0]}2^{[0]}$ & 0 & 1,2 &     &     & 1   &     &     & 1,2 &     &     & 1,2 &     &     \\ \hline
$1^{[0]}2^{[1]}$ & 4 & 1   &     &   2 & 1   &     &     & 1   &     &   2 & 1   &     &   2 \\ \hline
$1^{[0]}2^{[2]}$ & 3 & 1   &   2 &     & 1   &     &     & 1   &   2 &     & 1   &   2 &     \\ \hline
$1^{[1]}2^{[0]}$ & 2 &   2 &     & 1   &     &     & 1   &   2 &     & 1   &   2 &     & 1   \\ \hline
$1^{[1]}2^{[1]}$ & 6 &     &     & 1,2 &     &     & 1   &     &     & 1,2 &     &     & 1,2 \\ \hline
$1^{[1]}2^{[2]}$ & 5 &     &   2 & 1   &     &     & 1   &     &   2 & 1   &     &   2 & 1   \\ \hline
$1^{[2]}2^{[0]}$ & 1 &   2 & 1   &     &     & 1   &     &   2 & 1   &     &   2 & 1   &     \\ \hline
$1^{[2]}2^{[1]}$ & 5 &     & 1   &   2 &     & 1   &     &     & 1   &   2 &     & 1   &   2 \\ \hline
$1^{[2]}2^{[2]}$ & 4 &     & 1,2 &     &     & 1   &     &     & 1,2 &     &     & 1,2 &     \\ \hline \hline
$2^{[0]}1^{[0]}$ & 1 & 1   &     &     & 1,2 &     &     &   2 &     &     & 1   &     &     \\ \hline
$2^{[0]}1^{[1]}$ & 3 &     &     & 1   &   2 &     & 1   &   2 &     &     & 1   &     &     \\ \hline
$2^{[0]}1^{[2]}$ & 2 &     & 1   &     &   2 & 1   &     &   2 &     &     & 1   &     &     \\ \hline
$2^{[1]}1^{[0]}$ & 5 &     &     & 1   &     &     & 1,2 &     &     &   2 &     &     & 1   \\ \hline
$2^{[1]}1^{[1]}$ & 4 &     & 1   &     &     & 1   &   2 &     &     &   2 &     &     & 1   \\ \hline
$2^{[1]}1^{[2]}$ & 3 & 1   &     &     & 1   &     &   2 &     &     &   2 &     &     & 1   \\ \hline
$2^{[2]}1^{[0]}$ & 3 &     & 1   &     &     & 1,2 &   2 &     &   2 &     &     & 1   &     \\ \hline
$2^{[2]}1^{[1]}$ & 2 & 1   &     &     & 1   &   2 &     &     &   2 &     &     & 1   &     \\ \hline
$2^{[2]}1^{[2]}$ & 4 &     &     & 1   &     &   2 & 1   &     &   2 &     &     & 1   &     \\ \hline
\end{tabular}
\caption{$(\sor,\Cyc^t,\Lmic^t,\Lmal^t,\Lmap^t)$ for $\G_{3,2}$}
\end{table}
\end{tiny}
%
%
%

As for the generating function, if $\ff=(1,2)$, then $H(\ff)=(1,2)$.
Main Theorem B says that
\begin{align*}
&\sum_{\pi\in\G_{3,2,(1,2)}}q^{\ell(\pi)}\prod_{t=0}^{2}\left( \prod_{i\in\Rmil^{-t}(\pi)}x_{t,i}\prod_{i\in\Lmil^{-t}(\pi)}y_{t,i}\right)\\
=& \sum_{\pi\in\G_{3,2,(1,2)}}q^{\sor(\pi)}\prod_{t=0}^{2}\left(\prod_{i\in\Cyc^t(\pi)}x_{t,i}\prod_{i\in\Lmic^t(\pi)}y_{t,i}\right) \\
=& (x_{0,1}y_{0,1}+x_{2,1}y_{2,1}q+x_{1,1}y_{1,1}q^2) (x_{0,2}+x_{2,2}q^3+x_{1,2}q^4).
\end{align*}

In the above the $q^4$ term is $$1\cdot x_{0,1}x_{1,2}y_{0,1}q^4 + 1\cdot x_{2,1}x_{2,2}y_{2,1}q^4.$$
If we look at the second equality, this term says that there are two colored permutations $\pi\in\G_{3,2}$. 
The first permutation has $1\in\Cyc^0(\pi)$, $2\in\Cyc^1(\pi)$, $1\in\Lmic^0(\pi)$ and $\sor(\pi)=4$, namely the permutation $1^{[0]}2^{[1]}$, and the second one has $1,2\in \Cyc^2(\pi)$, $1\in \Lmic^2(\pi)$ and $\sor(\pi)=4$, namely the permutation $1^{[2]}2^{[2]}$. 
Similary from the first equality this term says that there are two colored permutations $\pi$: one has $1\in\Rmil^{0}(\pi)$, $2\in\Rmil^{-1}(\pi)=\Rmil^{2}(\pi)$, $1\in\Lmil^{0}(\pi)$ and $\ell(\pi)=4$, namely the permutation $1^{[0]}2^{[2]}$, and the other has $1,2\in\Rmil^{-2}(\pi)=\Rmil^{1}(\pi)$, $1\in\Lmil^{-2}(\pi)=\Lmil^{1}(\pi)$ and $\ell(\pi)=4$, namely the permutation $1^{[1]}2^{[1]}$.


\subsection{Even-signed permutations on a Ferrers shape}
The second part of the paper is to investigate Coxeter group of type D, or even-signed permutations. 
The sorting index $\sor_D$ was defined by~\cite{Petersen_11} and the set-valued equidistibution result restricted to a Ferrers shape was investigated in~\cite{Poznanovic_14}. 
Our two main results are Theorem~\ref{thm:equi_restricted_Dn} and Theorem~\ref{cor:GF_Dn}, which are the type $D$ version of Main Theorem A and Main Theorem B respectively, refining the results in~\cite{Poznanovic_14}.

\medskip
The rest of the paper is organized as follows. 
In Section~\ref{sec:def} we introduce the colored permutation groups $\G_{r,n}$ and define the sorting index and various set-valued Stirling statistics. 
A bijection on $\G_{r,n}$, again inspired by Foata and Han~\cite{Foata_09}, is given in Section~\ref{sec:biject}.
Based on this bijection, in Section~\ref{SectionMainA} and~\ref{SectionMainB} we prove Main Theorem A and B respectively.
In Section~\ref{sec:Dn} we investigate even-signed permutations.


\section{Sorting index of colored permutations}\label{sec:def}
\subsection{Colored permutations}

Let $r,n$ be positive integers. 
The group of \emph{colored permutations} $\G_{r,n}$ of $n$ letters with $r$ colors is $$\G_{r,n}:=C_r\wr\mathfrak{S}_n,$$ the wreath product of the cyclic group $C_r(:=\mathbb{Z}/r\mathbb{Z})$ with $\mathfrak{S}_n$.

An elements of $\G_{r,n}$ is an ordered pair $(\sigma,\mathbf{z})$, where $\sigma=\sigma_1\sigma_2\cdots\sigma_n\in\mathfrak{S}_n$ and $\mathbf{z}=(z_1,z_2,\ldots,z_n)$ is an $n$-tuple of integers with $z_i\in C_r$. 
The product of $(\sigma,\mathbf{z})$ and $(\rho,\mathbf{w})$ is $(\sigma\rho,\mathbf{w}+\rho(\mathbf{z}))$, where $\rho(\mathbf{z}):=(z_{\rho(1)},z_{\rho(2)},\ldots,z_{\rho(n)})$ and the addition is taken as mod $r$. 
It is easy to see that $e=(12\cdots n,(0,0,\ldots,0))$ is the identity.

We can represent elements of $\G_{r,n}$ in different ways. 
Let $$\Sigma:=\{1,\ldots,n,\bar{1},\ldots,\bar{n},\bar{\bar{1}},\ldots,\bar{\bar{n}},\ldots,1^{[r-1]},\ldots,n^{[r-1]}\},$$ 
then $(\sigma,\mathbf{z})$ can be viewed as the bijection $\pi:\Sigma\to\Sigma$ such that $\pi(i)=\sigma_i^{[z_i]}$ for $i\in [n]$ and $\pi(\bar{i})=\overline{\pi(i)}$ for $i\in\Sigma$.
For instance, $(3214,(2,1,1,0))\in\G_{3,4}$ can be represented as the bijection
$$
\left(
\begin{array}{cccccccccccc}
\bar{\bar{1}} & \bar{\bar{2}} & \bar{\bar{3}} & \bar{\bar{4}} & \bar{1} & \bar{2} & \bar{3} & \bar{4} & 1 & 2 & 3 & 4  \\
\bar{3} & 2 & 1 & \bar{\bar{4}} & 3 & \bar{\bar{2}} & \bar{\bar{1}} & \bar{4} & \bar{\bar{3}} & \bar{2} & \bar{1} & 4
\end{array}
\right),
$$
which is called its \emph{two-line notation}. 
By omitting the first row we have the \emph{one-line notation} $\bar{3}21\bar{\bar{4}}~3\bar{\bar{2}}\bar{\bar{1}}\bar{4}~\bar{\bar{3}}\bar{2}\bar{1}4$, or more tersely the \emph{window notation} $\bar{\bar{3}}\bar{2}\bar{1}4$ by only recording the image of $[n]$.
In this manner we write an element $(\sigma,\mathbf{z})\in\G_{r,n}$ as a word
$$\sigma_1^{[z_1]}\,\sigma_2^{[z_2]}\cdots\sigma_n^{[z_n]},$$
in which $\sigma_i$ and $z_i$ are respectively called the \emph{base value} and \emph{color} of $\pi(i)$.



The group $\G_{r,n}$ can be generated by the set of generators
$$\mathcal{S}_n=\{s_0,s_1,\ldots,s_{n-1}\}.$$ In the window notation,
$s_i$ ($1\le i\le n$) is the transpostition of swapping the $i$-th
and $(i+1)$-th letters, while $s_0$ is the action of adding one more
bar on the first letter (the number of bars is taken module $r$).
Note that the multiplication is on the right. For example, if
$r=3,n=4$, then
$$s_0s_1s_0s_2s_1s_0s_0s_3=\bar{\bar{3}}\bar{2}4\bar{1}.$$
The generators are subject to the conditions
\begin{equation} \label{eq:wreath_product}
\begin{cases}
s_0^r &= 1, \notag \\
s_i^2 &= 1, \text{ } 1\le i\le n-1, \notag \\
(s_is_j)^2 &= 1, \text{ } |i-j|>1, \notag \\
(s_is_{i+1})^3 &= 1, \text{ } 1\le i \le n-2, \\
(s_0s_1)^{2r} &= 1. \notag
\end{cases}
\end{equation}

\subsection{Length}\label{sec:def_mahonian}
We define the \emph{length} $\ell(\pi)$ of $\pi=\sigma_1^{[z_1]}\,\sigma_2^{[z_2]}\cdots\sigma_n^{[z_n]}\in\G_{r,n}$ to be the minimal number of generators in $\mathcal{S}_n$ needed to
represent it. 
Bagno~\cite{Bagno_04} gave the following combinatorial interpretation of $\ell(\pi)$:
\begin{equation}\label{eq:ell_defi}
\ell(\pi)=\inv(\pi)+\sum_{z_i>0}(\sigma_i+z_i-1),
\end{equation}
where $\inv(\pi)$ is the number of inversions in the window notation of $\pi$ with respect to the linear order
\begin{equation}\label{eq:order_wreath}
n^{[r-1]}<\cdots<\bar{n} <\cdots< 1^{[r-1]}<\cdots < \bar{1} < 1<\cdots<n.
\end{equation}
For example, let $\pi=\bar{\bar{3}}\bar{2}4\bar{1}\in\G_{3,4}$. 
Then $\ell(\pi)=1+(4+2+1)=8$. 
The distribution of $\ell$ was also derived in \cite{Bagno_04} as
\begin{equation}\label{eq:distribution_ell}
\sum_{\pi\in\G_{r,n}}q^{\ell(\pi)} = [n]_q!\cdot \prod_{i=1}^n(1+q^i[r-1]_q).
\end{equation}

\medskip

There is another `length function' on $\G_{r,n}$. 
Consider the one-line notation. 
For $1\leq i<j\leq n$ and $0\leq t<r$ let $(i^{[t]}\,j)$ be the transposition of swapping the $i^{[t]}$-th with $j$-th letters, the $i^{[t+1]}$-th with $\bar{j}$-th letters, $\ldots$, the $i^{[t+r-1]}$-th with $j^{[r-1]}$-th letters.
Also, for $1\leq i\leq n$ and $0<t<r$ let $(i^{[t]}\,i)$ be the action of adding $t$ bars on the $i$-th, $\bar{i}$-th, $\ldots$, $i^{[r-1]}$-th letters.
In the window notation, multiplying $\pi=\sigma_1^{[z_1]}\,\sigma_2^{[z_2]}\cdots\sigma_n^{[z_n]}\in\G_{r,n}$ on the right by $(i^{[t]}\,j)$, $i<j$, has the effect of replacing $\pi_j$ by $\sigma_i^{[z_i+t]}$ and $\pi_i$ by $\sigma_j^{[z_j-t]}$, while multiplying $\pi$ on the right by $(i^{[t]}\,i)$ has the effect of replacing $\pi_i$ by $\sigma_i^{[z_i+t]}$.
For example, if $\pi = 2\bar{\bar{5}}1\bar{4}\bar{\bar{3}}\in\G_{3,5}$, then $\pi\cdot(\bar{2}\,5) = 2\bar{3}1\bar{4}5$ and $\pi\cdot(\bar{\bar{5}}\,5) = 2\bar{\bar{5}}1\bar{4}\bar{3}$.

It can be seen that $\G_{r,n}$ can also be generated by
$$\mathcal{T}_n:=\{(i^{[t]}\,j):\,1\leq i<j\leq n \text{ and } 0\leq t<r\} \cup \{(i^{[t]}\,i):\,1\leq i\leq n \text{ and } 0< t<r\}.$$

Denote by $\ell'(\pi)$ the minimal number of elements in $\mathcal{T}_n$ needed to express $\pi$. 
The distribution of $\ell'$ will be derived in Corollary~\ref{cor:distribution_ell'} as
$$\sum_{\pi\in\G_{r,n}}t^{\ell'(\pi)} = \prod_{i=1}^n(1+(ri-1)t).$$

Note that $\ell'$ is also called the \emph{reflection length} when $r=1,2$~\cite{Petersen_11}, where each element of $\mathcal{T}_n$ is a reflection. However when $r\ge 3$ elements of $\mathcal{T}_n$ are not reflections.

\subsection{Sorting index}
Now we come to the key definition of the whole paper. We will define
a reasonable sorting index on $\G_{r,n}$. Note that any
$\pi\in\G_{r,n}$ can be uniquely written as a product
$$\pi=(i_1^{[t_1]}\,j_1)(i_2^{[t_2]}\,j_2)\cdots (i_k^{[t_k]}\,j_k)$$ for
some $k$ such that $0<j_1<\cdots < j_k$.
\begin{defi} The \emph{sorting index} of $\pi\in\G_{r,n}$ is
\begin{equation}\label{eq:sor_defi}
\sor(\pi) = \sum_{s=1}^k \Big(j_s - i_s + \chi(t_s>0)\cdot \big(2(i_s-1)+t_s\big)\Big),
\end{equation}
where $\chi(\text{A})=1$ if the statement $\text{A}$ is true, or
$\chi(\text{A})=0$ otherwise.
\end{defi}

It can be checked that when $r=1,2$ our definitions meet the
definitions in type A and B in~\cite{Petersen_11}.

$\sor(\pi)$ can be computed conveniently on the labeled graph $G^{(r,n)}$ in the way we explain below. 
The graph $G^{(r,n)}$ has $r\times n$ vertices arranged in an rectangle shape. 
Two vertices are connected by an edge if they are adjacent and are of the same row or on the leftmost column.
For convenience, the columns are indexed by $1,2,\ldots,n$ and the rows by $0,1,\ldots,r-1$.

\begin{figure}[h]
\includegraphics[width=3in]{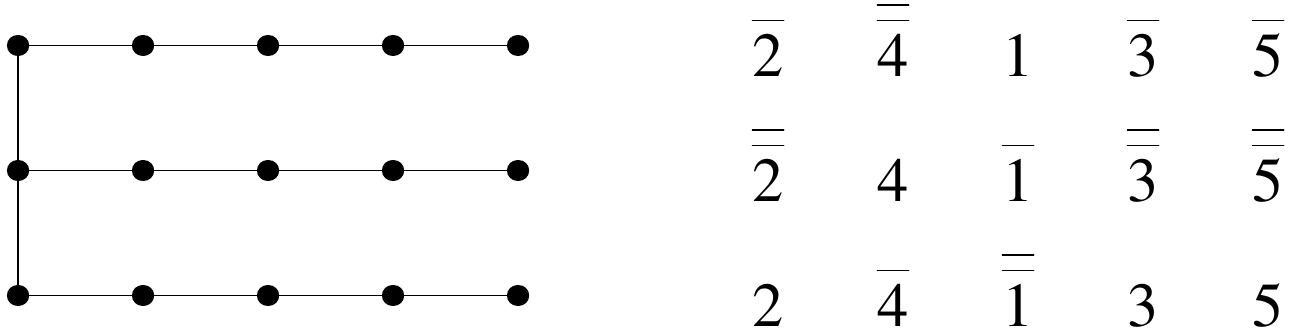}
\caption{The $G^{(3,5)}$ with the labeling with respect to $\bar{2}\bar{\bar{4}}1\bar{3}\bar{5}$.} \label{fig:comb}
\end{figure}

For $1\le i\le n$ and $0\le j\le r-1$, the vertice of the $i$-th column and $j$-row is labelled by $\pi(i^{[j]})$. 
Simply put, we label vertices of the first row by $\pi$, and the labels of the $j$-th row are obtained by adding a bar on each letter of the $(j-1)$-th row. 
See Fig.~\ref{fig:comb} for the graph $G^{(3,5)}$ and its labelling from $\pi=\bar{2}\bar{\bar{4}}1\bar{3}\bar{5}$.

The sorting process of $\pi$ can be done on $G^{(r,n)}$ in the following way. 
The goal of the sorting is the identity permutation $(12\cdots n,(0,0,\ldots,0))$, or $12\cdots n$ for short.\\
\begin{enumerate}
\item[(i)] Find the largest unsorted letter $j$ in $\{1,2,\dots, n\}$. Suppose it is at the $i$-th column and the $t$-th row.
\item[(ii)] Exchange $j$ with the $j$-th letter of the $0$-th row. The distance for $j$ needed to travel along the graph $G^{(r,n)}$ to its new position (namely, row $0$, column $j$) is recorded.
\item[(iii)] If $i\ne j$, relabel the vertices of the $i$-th and $j$-th columns by fixing the two exchanged letters first and then following the labelling rule above. If $i=j$, relabel the $j$-th column by fixing $j$ first and then following the labelling rule.
\item[(iv)] Back to (i).
\end{enumerate}
Then the sorting index $\sor(\pi)$ is the total distances in the above process. 
For example, the process for sorting $\pi=\bar{2}\bar{\bar{4}}1\bar{3}\bar{5}\in\G_{3,5}$ can be illustrated as

$$
\begin{array}{cccccc}
 &
    \begin{array}{ccccc}
    \bar{2} & \bar{\bar{4}} & 1 & \bar{3} & \bar{5} \\
    \bar{\bar{2}} & 4 & \bar{1} & \bar{\bar{3}} & \bar{\bar{5}} \\
    2 & \bar{4} & \bar{\bar{1}} & 3 & {\bf 5}
    \end{array}
 &
    \stackrel{(\bar{\bar{5}}\,5)}{\longrightarrow}
 &
    \begin{array}{ccccc}
    \bar{2} & \bar{\bar{4}} & 1 & \bar{3} & 5 \\
    \bar{\bar{2}} & {\bf 4} & \bar{1} & \bar{\bar{3}} & \bar{5} \\
    2 & \bar{4} & \bar{\bar{1}} & 3 & \bar{\bar{5}}
    \end{array}
 &
    \stackrel{(\bar{2}\,4)}{\longrightarrow}
 &
    \begin{array}{ccccc}
    \bar{2} & {\bf 3} & 1 & 4 & 5 \\
    \bar{\bar{2}} & \bar{3} & \bar{1} & \bar{4} & \bar{5} \\
    2 & \bar{\bar{3}} & \bar{\bar{1}} & \bar{\bar{4}} & \bar{\bar{5}}
    \end{array}
\\ & & & & & \\
    \stackrel{(2\,3)}{\longrightarrow}
 &
    \begin{array}{ccccc}
    \bar{2} & 1 & 3 & 4 & 5 \\
    \bar{\bar{2}} & \bar{1} & \bar{3} & \bar{4} & \bar{5} \\
    {\bf 2} & \bar{\bar{1}} & \bar{\bar{3}} & \bar{\bar{4}} & \bar{\bar{5}}
    \end{array}
 &
    \stackrel{(\bar{\bar{1}}\,2)}{\longrightarrow}
 &
    \begin{array}{ccccc}
    \bar{1} & 2 & 3 & 4 & 5 \\
    \bar{\bar{1}} & \bar{2} & \bar{3} & \bar{4} & \bar{5} \\
    {\bf 1} & \bar{\bar{2}} & \bar{\bar{3}} & \bar{\bar{4}} & \bar{\bar{5}}
    \end{array}
 &
    \stackrel{(\bar{1}\,1)}{\longrightarrow}
 &
    \begin{array}{ccccc}
    1 & 2 & 3 & 4 & 5 \\
    \bar{1} & \bar{2} & \bar{3} & \bar{4} & \bar{5} \\
    \bar{\bar{1}} & \bar{\bar{2}} & \bar{\bar{3}} & \bar{\bar{4}} & \bar{\bar{5}}
    \end{array}
\end{array}.
$$
Therefore,
$\pi=(\bar{1}\,1)(\bar{\bar{1}}\,2)(2\,3)(\bar{2}\,4)(\bar{\bar{5}}\,5)$
and $\sor(\pi)=2+3+1+5+10=21$.

It is easy to see that to apply (i),(ii),(iii) once is equivalent to
multiply the permutation by the transposition $(i^{[t]}\,j)$.

\subsection{Sorting index v.s. length}
There is a simple way to generate elements of $\G_{r,n}$
recursively. The idea is to append the letter $n$ to the end of
$\pi=\pi_1\cdots\pi_{n-1}\in\G_{r,n-1}$, then pick a letter $i$ and
a color $t$ arbitrary and apply the transposition $(i^{[t]}\,n)$.
More formally, we define elements $\Phi_i$ of the group algebra of
$\G_{r,n}$ by $\Phi_1:=1+(\bar{1}\,1)+\cdots +(1^{[r-1]}\,1)$ and
for $2\le j\le n$
$$\Phi_j:=1+\sum_{i=1}^{j-1}(i\,j)+\sum_{t=1}^{r-1}\sum_{i=1}^j(i^{[t]}\,j).$$

\begin{lem}\label{lem:factorization}
We have
$$\Phi_1\Phi_2\cdots\Phi_n=\sum_{\pi\in\G_{r,n}}\pi.$$
\end{lem}
\proof Obviously the formula holds for $n=1$. Suppose by induction
that $\Phi_1\Phi_2\cdots\Phi_{n-1}=\sum_{\pi\in\G_{r,n-1}}\pi$.
Observe that $\Phi_n$ is $1$ plus the sum of all transpositions in
$\mathcal{T}$ involving the letter $n$. Thus, for
$\pi=\pi_1\cdots\pi_{n-1}n\in\G_{r,n}$ in the window notation we have
$$
\begin{array}{rllll}
\pi\cdot\Phi_n & =\pi_1\cdots\pi_{n-1}n & +~\pi_1\cdots n\pi_{n-1} & +~\cdots & +~n\pi_2\cdots\pi_{n-1}\pi_1 \\ \\
 & +~\pi_1\cdots\pi_{n-1}n^{[r-1]} & +~\pi_1\cdots n^{[r-1]}\pi_{n-1}^{[1]} & +~\cdots & +~n^{[r-1]}\pi_2\cdots\pi_{n-1}\pi_1^{[1]} \\ \\
 & +~\cdots\cdots\cdots & & & \\ \\
 & +~\pi_1\cdots\pi_{n-1}n^{[1]} & +~\pi_1\cdots n^{[1]}\pi_{n-1}^{[r-1]} & +~\cdots & +~n^{[1]}\pi_2\cdots\pi_{n-1}\pi_1^{[r-1]}.
\end{array}
$$
It is clear now that $\pi\cdot\Phi_n=\pi'\cdot\Phi_n$ iff $\pi=\pi'$
and hence we have
$$\Phi_1\Phi_2\cdots\Phi_{n-1}\Phi_n = \sum_{\pi\in\G_{r,n},\,\pi(n)=n}\pi\cdot\Phi_n = \sum_{\pi\in\G_{r,n}}\pi.$$
\qed

Now we prove the main result of this section.

\begin{thm} \label{thm:sor_GF}
The statistics $\ell$ and $\sor$ have the same distribution over $\G_{r,n}$. That is,
\begin{eqnarray*}
\sum_{\pi\in\G_{r,n}}q^{\ell(\pi)} &=& \sum_{\pi\in\G_{r,n}}q^{\sor(\pi)} \\
&=& [n]_q!\cdot \prod_{i=1}^n(1+q^i[r-1]_q).
\end{eqnarray*}
\end{thm}

\proof We proceed by induction. Define the linear mapping
${\it\Phi}:\mathbb{Z}(\G_{r,n})\to\mathbb{Z}(q)$ by
$${\it\Phi}(\pi):=q^{\sor(\pi)}.$$ It is obvious that
\begin{align*}
{\it\Phi}(\Phi_i) =& (1+q+\cdots+q^{i-1})+(q^i+q^{i+1}+\cdots+q^{2i-1})+(q^{i+1}+q^{i+2}+\cdots+q^{2i})\\
&+ \cdots + (q^{i+r-2}+q^{i+r-1}+\cdots+q^{2i+r-3})\\
=& \left(1+q^i[r-1]_q\right) [i]_q.
\end{align*}
Thus, by Lemma~\ref{lem:factorization}, it suffices to show that
$${\it\Phi}(\Phi_1\cdots\Phi_{n-1}){\it\Phi}(\Phi_n) = {\it\Phi}(\Phi_1\cdots\Phi_{n-1}\Phi_n).$$
Following the proof of Lemma~\ref{lem:factorization}, let
$\pi=\pi_1\cdots\pi_{n-1}n\in\G_{r,n}$ be in the window notation. Since
$\pi(n)=n$, by the definition of $\sor$ we have $\sor(\pi\cdot
(i\,n))=\sor(\pi)+(n-i)$ for $1\leq i<n$, and $\sor(\pi\cdot
(i^{[t]}\,n))=\sor(\pi)+(n+t-i-2)$ for $1\leq i\leq n$ and $1\leq
t<r$. This implies that
\begin{align*}
{\it\Phi}(\pi\cdot\Phi_n) =& {\it\Phi}(\pi) \Big(1 + \sum_{i=1}^{n-1}q^{n-i} + \sum_{t=1}^{r-1}\sum_{i=1}^n q^{n+t+i-2}\Big) \\
=& {\it\Phi}(\pi) \Big(\sum_{i=0}^{n-1} q^i + q^n \sum_{t=1}^{r-1}q^{t-1} \sum_{i=1}^{n}q^{i-1}\Big) \\
=& {\it\Phi}(\pi) \Big(1+q^n[r-1]_q\Big)[n]_q \\
=& {\it\Phi}(\pi){\it\Phi}(\Phi_n).
\end{align*}
Thus,
\begin{align*}
{\it\Phi}(\Phi_1\cdots\Phi_{n-1}\Phi_n) =& {\it\Phi}\Big(\sum_{\pi\in\G_{r,n},\,\pi(n)=n}\pi\cdot\Phi_n\Big)\\
=& \sum_{\pi\in\G_{r,n},\,\pi(n)=n} {\it\Phi}(\pi\cdot\Phi_n) \\
=& {\it\Phi}(\Phi_n) \sum_{\pi\in\G_{r,n},\,\pi(n)=n} {\it\Phi}(\pi)\\
=& {\it\Phi}(\Phi_n) {\it\Phi}(\Phi_1\cdots\Phi_{n-1}),
\end{align*}
as desired.
\qed

\subsection{Set-valued Stirling statistics}\label{sec:def_stirling}
We can also represent a colored permutation
$\pi=\sigma_1^{[z_1]}\sigma_2^{[z_2]}\cdots\sigma_n^{[z_n]}\in
\G_{r,n}$ as a product of disjoint colored cycles. This is done by
write $\sigma_1\sigma_2\cdots\sigma_n$ (taken as a permutation
in $\mathfrak{S}_n$) in its cycle decomposition and add back bars
of each letters. For example, in $G_{3,9}$, we have
$$\bar{5}\bar{\bar{6}}\bar{3}\bar{1}4\bar{\bar{2}}79\bar{\bar{8}} =
(\bar{1}\,\bar{5}\,4)(\bar{\bar{2}}\bar{\bar{6}})(\bar{3})(7)(\bar{\bar{8}}\,9).$$
We note that this is \emph{not} the same with the cycle
decomposition from the two-line notation.

Assume that $C_1,C_2,\ldots,C_k$ are the disjoint colored cycles of
$\pi$. For
$C_i=(\sigma_{i_1}^{[z_{i_1}]}\sigma_{i_2}^{[z_{i_2}]}\cdots)$, let
$\alpha_i$ be the smallest integer in
$\{\sigma_{i_1},\sigma_{i_2},\ldots\}$ and $c_i$ be the remainder of
$\sum z_{i_j}$ divided by $r$. We define the \emph{colored cycle
set} of $\pi$ by
$$\Cyc(\pi):=\{\alpha_1^{[c_1]},\ldots,\alpha_k^{[c_k]}\},$$
and for each $t=0,1,\ldots,r-1$ the \emph{refined colored cycle
set} by
$$\Cyc^t(\pi):=\{\alpha_i:\,\alpha_i^{[c_i]}\in\Cyc(\pi)
\text{ and } c_i=t\}.$$ Also we denote by $\cyc$ and $\cyc^t$,
$0\leq t\leq r-1$, the cardinalities of $\Cyc$ and $\Cyc^t$
respectively. In the above example,
$\pi=\bar{5}\bar{\bar{6}}\bar{3}\bar{1}4\bar{\bar{2}}79\bar{\bar{8}}\in\G_{3,9}$
has $\Cyc(\pi)=\{\bar{\bar{1}},\bar{2},\bar{3},7,\bar{\bar{8}}\}$,
$\Cyc^0(\pi)=\{7\}$, $\Cyc^1(\pi)=\{2,3\}$ and
$\Cyc^2(\pi)=\{1,8\}$.

In Corollary~\ref{cor:equidistribution_stirling} we will prove that
\begin{equation}
\cyc^0(\pi)=n-\ell'(\pi),
\end{equation}
and in Corollary~\ref{cor:distribution_ell'} the distribution of
$\cyc^0$ will be derived as
\begin{equation}
\sum_{\pi\in\G_{r,n}}t^{\cyc^0(\pi)} = \prod_{i=1}^n(t+ri-1),
\end{equation}
which are counterparts of (\ref{eq:cyc_ell'}) and
(\ref{cycle_distribution}).
%
We may call a statistic \emph{colored Stirling} if it is
equidistributed with $\cyc^0$ over $\G_{r,n}$.

Now we define the following set-valued statistics for $\pi=\sigma_1^{[z_1]}\sigma_2^{[z_2]}\cdots\sigma_n^{[z_n]}\in\G_{r,n}$ :
\begin{enumerate}
\item $\Rmil$, the set of \emph{right-to-left minimum letters}:
    $$\Rmil(\pi):=\{\sigma_i^{[z_i]}:\, \sigma_i < \sigma_j \text{ for any } j > i\}.$$
\item $\Rmip$, the set of \emph{right-to-left minimum places}:
    $$\Rmip(\pi):=\{i^{[z_i]}:\, \sigma_i < \sigma_j \text{ for any } j > i\}.$$
\item $\Lmal$, the set of \emph{left-to-right maximum letters}:
    $$\Lmal(\pi):=\{\sigma_i^{[z_i]}:\, \sigma_i > \sigma_j \text{ for any } j < i\}.$$
\item $\Lmap$, the set of \emph{left-to-right maximum places}:
    $$\Lmal(\pi):=\{i^{[z_i]}:\, \sigma_i > \sigma_j \text{ for any } j < i\}.$$
\item $\Lmil$, the set of \emph{left-to-right minimum letters}:
    $$\Lmil(\pi):=\{\sigma_i^{[z_i]}:\, \sigma_i < \sigma_j \text{ for any } j < i\}.$$
\item $\Lmic$, the set of \emph{left-to-right minimum letters in the first cycle}:
    $$\Lmic(\pi):=\Lmil\big(\pi(1)\,\pi^2(1)\,\pi^3(1)\cdots\big).$$
\end{enumerate}
Let $\rmin,\lmin,\lmax$ and $\lmic$ denote the cardinalities of
$\Rmil$, $\Lmil$, $\Lmal (\mbox{or }\Lmap)$ and $\Lmic$
respectively. For each $t=0,1,\ldots, r-1$ the statistics
$\Rmil^t,\Lmil^t,\Lmal^t,\Lmap^t,\Lmic^t$ and
$\rmin^t,\lmin^t,\lmax^t,\lmic^t$ are defined similarly. In
Corollary~\ref{cor:equidistribution_stirling} we will prove that
$\cyc^0,\rmin^0,\lmin^0,\lmax^0$ and $\lmic^0$ are colored Stirling
and $\cyc,\rmin,\lmin,\lmax$ and $\lmic$ are equidistributed over
$\G_{r,n}$.

For example, Table 3 lists these set-valued statistics for
$\pi=\bar{5}\bar{\bar{6}}\bar{3}\bar{1}4\bar{\bar{2}}79\bar{\bar{8}}\in\G_{3,9}$.
Note that $\pi(1)\,\pi^2(1)\,\pi^3(1)\cdots =
\bar{5}\bar{4}\bar{\bar{1}}54\bar{1}\bar{\bar{5}}\bar{\bar{4}}1$.

\footnotesize
\begin{table}[h]
\begin{tabular} {c||c|c|c|c|c|c|c}
\hline
$(\Stat)$ & $\Cyc$ & $\Rmil$ & $\Rmip$ & $\Lmal$ & $\Lmap$ & $\Lmil$ & $\Lmic$ \\
\hline \hline \rule{0pt}{10pt}
$\Stat(\pi)$ & $\bar{\bar{1}},\bar{2},\bar{3},7,\bar{\bar{8}}$ & $\bar{1},\bar{\bar{2}},7,\bar{\bar{8}}$ & $\bar{4},\bar{\bar{6}},7,\bar{\bar{9}}$ & $\bar{5},\bar{\bar{6}},7,\bar{\bar{8}}$ & $\bar{1},\bar{\bar{2}},7,\bar{\bar{9}}$ & $\bar{1},\bar{3},\bar{5}$ & $\bar{\bar{1}},\bar{4},\bar{5}$ \\
\hline \rule{0pt}{10pt}
$\Stat^0(\pi)$ & $7$ & $7$ & $7$ & $7$ & $7$ & $\emptyset$ & $\emptyset$ \\
\hline \rule{0pt}{10pt}
$\Stat^1(\pi)$ & $2,3$ & $1$ & $4$ & $5$ & $1$ & $1,3,5$ & $4,5$ \\
\hline \rule{0pt}{10pt}
$\Stat^2(\pi)$ & $1,8$ & $2,8$ & $6,9$  & $6,8$ & $2,9$ & $\emptyset$ & $1$ \\
\hline
\end{tabular}
\caption{Set-valued statistics for
$\pi=\bar{5}\bar{\bar{6}}\bar{3}\bar{1}4\bar{\bar{2}}79\bar{\bar{8}}$.}
\end{table} \label{tal:example}
\normalsize

\section{A bijection on $G_{r,n}$ }\label{sec:biject}

In this section we establish a bijection $\phi$ on $\G_{r,n}$, which
is the key ingredient for proving the Main Theorem A. It turns out
that $\phi$ is the composition
$$\phi:=(\text{B-code})^{-1}\circ \text{(A-code)}$$
of the A-code and B-code defined below. Note that $\phi$ is a
generalization of the bijection first defined by Foata and
Han~\cite{Foata_09} and extended by Chen-Guo-Gong~\cite{Chen_13} and
Poznanovi\'{c}~\cite{Poznanovic_14}.

\subsection{The A-code}\label{sec:biject_Acode}

For $\pi=\sigma_1^{[z_1]}\sigma_2^{[z_2]}\cdots\sigma_n^{[z_n]}\in\G_{r,n}$ define its \emph{Lehmer code} to be the sequence
$$\Leh(\pi)=(h_1^{[-z_1]},h_2^{[-z_2]},\ldots,h_n^{[-z_n]}),$$
where $-z_i$ is taken modulo $r$ and for each $i$
$$h_i:=|\{j:\,1\leq j\leq i \text{ and } \sigma_j\leq \sigma_i\}|.$$
And the \emph{A-code} of $\pi$ is then defined as
$$\text{A-code}(\pi):=\Leh(\pi^{-1}).$$

For example, for
$\pi=\bar{5}\bar{\bar{6}}\bar{3}\bar{1}4\bar{\bar{2}}79\bar{\bar{8}}\in\G_{3,9}$,
we have
$\pi^{-1}=\bar{\bar{4}}\bar{6}\bar{\bar{3}}5\bar{\bar{1}}\bar{2}7\bar{9}8$
and $\text{A-code}(\pi)=
(\bar{1},\bar{\bar{2}},\bar{1},3,\bar{1},\bar{\bar{2}},7,\bar{\bar{8}},8)$.
It is obvious that the A-code is a bijection from $\G_{r,n}$ to the
set
$$\mathsf{CS}_{r,n}:=\{(c_1^{[e_1]},c_2^{[e_2]},\ldots,c_n^{[e_n]})
: 1\leq c_i\leq i, 0 \leq e_i<r \mbox{ for each }
1\le i\le n\}.$$

\medskip
We can compute A-code by the following algorithm.\\
 \noindent {\bf Algorithm for A-code.} For $\pi$ we construct a
sequence of $n$ colored permutations
$\pi=\pi^{(n)},\pi^{(n-1)},\ldots,\pi^{(1)}$ such that
$\pi^{(i)}\in\G_{r,i}$ and meanwhile build up the A-code
$(a_1,a_2,\dots ,a_n)$, where $a_i$ is of the form $c_i^{[e_i]}$.
For $i$ from $n$ down to $2$ we look at $i^{[t]}$. If $i^{[t]}$
appears at the $p$-th position in $\pi^{(i)}$, then we set
$a_i=p^{[t]}$ and let $\pi^{(i-1)}$ be obtained from $\pi^{(i)}$ by
deleting the element $i^{[t]}$. Finally we let $a_1:=\pi^{(1)}(1)$.
It is easy to see the procedure is reversible.

For instance, for
$\pi=\bar{5}\bar{\bar{6}}\bar{3}\bar{1}4\bar{\bar{2}}79\bar{\bar{8}}\in\G_{3,9}$
we successively get
$$\begin{array}{llll}
\pi^{(9)} = \bar{5}\bar{\bar{6}}\bar{3}\bar{1}4\bar{\bar{2}}7{\bf 9}\bar{\bar{8}}, & p=8, & t=0, & a_9=8;\\
\pi^{(8)} = \bar{5}\bar{\bar{6}}\bar{3}\bar{1}4\bar{\bar{2}}7{\bf\bar{\bar{8}}}, & p=8, & t=2, & a_8=\bar{\bar{8}};\\
\pi^{(7)} = \bar{5}\bar{\bar{6}}\bar{3}\bar{1}4\bar{\bar{2}}{\bf 7}, & p=7, & t=0, & a_7=7;\\
\pi^{(6)} = \bar{5}{\bf\bar{\bar{6}}}\bar{3}\bar{1}4\bar{\bar{2}}, & p=2, & t=2, & a_6=\bar{\bar{2}};\\
\pi^{(5)} = {\bf\bar{5}}\bar{3}\bar{1}4\bar{\bar{2}}, & p=1, & t=1, & a_5=\bar{1};\\
\pi^{(4)} = \bar{3}\bar{1}{\bf 4}\bar{\bar{2}}, & p=3, & t=0, & a_4=3;\\
\pi^{(3)} = {\bf\bar{3}}\bar{1}\bar{\bar{2}}, & p=1, & t=1, & a_3=\bar{1};\\
\pi^{(2)} = \bar{1}{\bf\bar{\bar{2}}}, & p=2, & t=2, & a_2=\bar{\bar{2}};\\
\pi^{(1)} = {\bf\bar{1}}, & p=1, & t=1, & a_1=\bar{1}.\\
\end{array}$$

Thus
$\text{A-code}(\pi)=(\bar{1},\bar{\bar{2}},\bar{1},3,\bar{1},\bar{\bar{2}},7,\bar{\bar{8}},8)$.

We can read the length from A-code.
\begin{lem}\label{lem:Acode_ell}
Suppose
$\text{A-code}(\pi)=(c_1^{[e_1]},c_2^{[e_2]},\ldots,c_n^{[e_n]})$.
Then we have
$$\ell(\pi)=\sum_{i=1}^n \Big(i-c_i+\chi(e_i>0)\cdot\big(2(c_i-1)+e_i\big)\Big).$$
\end{lem}
\proof Consider the procedure of recovering the colored permutation
from its A-code. At the $i$-th step we insert the entry $i^{[e_i]}$
into the $c_i$-th position of $\pi^{(i-1)}$. From the definition of
$\ell$, it can be seen that after the $i$-th step the length
function increases by $i-c_i$ when $e_i=0$ and by $c_i-1+(i+e_i-1)$
when $e_i>0$. Hence we have
$$\ell(\pi^{(i)})-\ell(\pi^{(i-1)}) = i-c_i
+\chi(e_i>0)\cdot(2c_i-2+e_i).$$ Since $\ell(\pi^{(0)})=0$, the
result follows. \qed

\medskip

For
$a=(c_1^{[e_1]},c_2^{[e_2]},\ldots,c_n^{[e_n]})\in\mathsf{CS}_{r,n}$,
we define the set-valued statistics
$$\Max(a):=\{i^{[e_i]}:\,c_i=i\},\qquad  \Min(a):=\{i^{[e_i]}:\,c_i=1\},$$
and their refined versions
$$\Max^t(a):=\{i:\,c_i=i,e_i=t\}, \qquad \Min^t(a):=\{i:\,c_i=1,e_i=t\}$$
for each $0\le t \le r-1$.  $\Rmil(a)$, $\Rmip(a)$, $\Rmil^t(a)$ and
$\Rmip^t(a)$ are defined similarly by regarding $a$ as a word.

\begin{lem}\label{lem:Acode_stirling}
Let $\pi\in\G_{r,n}$ and $a=\text{A-code}(\pi)$. Then for each $0\le
t\le r-1$ we have
\begin{enumerate}
\item $\Rmil(\pi)=\Max(a)$ and $\Rmil^t(\pi)=\Max^{t}(a)$,
\item $\Lmil(\pi)=\Min(a)$ and $\Lmil^t(\pi)=\Min^{t}(a)$,
\item $\Lmap(\pi)=\Rmil(a)$ and $\Lmap^t(\pi)=\Rmil^{t}(a)$,
\item $\Lmal(\pi)=\Rmip(a)$ and $\Lmal^t(\pi)=\Rmip^{t}(a)$.
\end{enumerate}
\end{lem}
\proof
Let $\pi=(\sigma,\mathbf{z})=\sigma_1^{[z_1]}\cdots\sigma_n^{[z_n]}$ and $a=(c_1^{[e_1]},c_2^{[e_2]},\ldots,c_n^{[e_n]})$.

(1) Following the algorithmic construction of A-code,
$\sigma_i^{[z_i]}$ is a right-to-left minimum letter iff it is at
the last position in $\pi^{(\sigma_i)}$, which implies
$c_{\sigma_i}=\sigma_i$ and $e_{\sigma_i}=z_i$. Hence
$\Rmil(\pi)=\Max(a)$.

(2) Similar to (1), now $\sigma_i^{[z_i]}$ is a left-to-right
minimum letter iff it is at the first position in
$\pi^{(\sigma_i)}$, which implies $c_{\sigma_i}=1$ and
$e_{\sigma_i}=z_i$. Hence $\Lmil(\pi)=\Min(a)$.

(3) Suppose $\pi^{-1}=(\rho,\mathbf{w})$. Since
$(\sigma\rho,\mathbf{w}+\rho(\mathbf{z}))=\pi\circ\pi^{-1}=(12\cdots
n,(0,0,\ldots,0))$, we obtain $\rho=\sigma^{-1}$ and
$w_i=-z_{\rho(i)}$ for all $i$. In other words,
$\pi^{-1}=(\sigma^{-1},-\sigma^{-1}(\mathbf{z}))$, where
\begin{equation}\label{eq:pi_inverse}
-\sigma^{-1}(\mathbf{z}) = \big(-z_{\sigma^{-1}(1)},-z_{\sigma^{-1}(2)},\ldots,-z_{\sigma^{-1}(n)}\big).
\end{equation}
Since $\Lmap(\sigma)=\Rmil(\sigma^{-1})$, this implies that
$i^{[z_i]}$ is an element of $\Lmap(\pi)$ iff $i^{[t]}$ is an
element of $\Rmil(\pi^{-1})$ for some $t$, hence by
\eqref{eq:pi_inverse} we have $t=-z_{\sigma^{-1}(\sigma(i))}=-z_i$.
Therefore,
\begin{equation}\label{eq:pi_inverse_1}
i^{[z_i]}\in\Lmap(\pi) \mbox{ iff } i^{[-z_i]}\in\Rmil(\pi^{-1}).
\end{equation}
On the other hand by the definition of Lehmer code it is obvious
that
\begin{equation}\label{eq:pi_inverse_2}
\sigma_i^{[z_i]}\in\Rmil(\pi) \mbox{ iff } \sigma_i^{[-z_i]}\in\Rmil(\Leh(\pi)).
\end{equation}
Combining \eqref{eq:pi_inverse_1} and \eqref{eq:pi_inverse_2} we obtain
$$i^{[z_i]}\in\Lmap(\pi) \mbox{ iff } i^{[z_i]}\in\Rmil(\Leh(\pi^{-1}))=\Rmil(a),$$
hence $\Lmap(\pi)=\Rmil(a)$.

(4) Similar to the proof of (3). \qed

\medskip


\subsection{The B-code}\label{sec:biject_Bcode}
The \emph{B-code} of $\pi\in\G_{r,n}$ is defined in the following
way. For $i=1,2,,\ldots,n$ let $k_i$ be the smallest integer $k\geq
1$ such that the base value of $\pi^{-k}(i)$ is less than or equal
to $i$. Then we define
$$\text{B-code}(\pi):=(b_1,b_2,\ldots,b_n) \text{ with }
b_i=\pi^{-k_i}(i).$$

It is not difficult to see that B-code is a bijection from
$\G_{r,n}$ to $\mathsf{CS}_{r,n}$. A simple way to compute
$\text{B-code}(\pi)$ is from its cycle representation. For example,
the cycle representation of
$\pi=\bar{5}\bar{\bar{6}}\bar{3}\bar{1}4\bar{\bar{2}}79\bar{\bar{8}}\in\G_{3,9}$
is
$$(1\,\bar{5}\,\bar{4}\,\bar{\bar{1}}\,5\,4\,\bar{1}\,\bar{\bar{5}}\,\bar{\bar{4}})(2\,\bar{\bar{6}}\,\bar{2}\,6\,\bar{\bar{2}}\,\bar{6})(3\,\bar{3}\,\bar{\bar{3}})(7)(\bar{7})(\bar{\bar{7}})(8\,9\,\bar{\bar{8}}\,\bar{\bar{9}}\,\bar{8}\,\bar{9})$$
and
$\text{B-code}(\pi)=(\bar{1},\bar{\bar{2}},\bar{\bar{3}},\bar{\bar{1}},\bar{\bar{1}},\bar{2},7,\bar{8},8)$.

\medskip
We can also compute B-code by the following algorithm.

\noindent {\bf Algorithm for B-code.} For $\pi\in\G_{r,n}$ we
construct a sequence of colored permutations
$\pi=\pi^{(n)},\pi^{(n-1)},\ldots,\pi^{(1)}$ such that
$\pi^{(i)}\in\G_{r,i}$ and meanwhile build up the B-code
$(b_1,b_2,\ldots,b_n)$. For $i$ from $n$ down to $2$ we assume that
$\pi^{(i)}(p^{[t]})=i$ for some $p$ and $t$. Set $b_i=p^{[t]}$ and
$\pi'=\pi^{(i)}\cdot(p^{[t]}\,i)$, the product of $\pi^{(i)}$ and
the transposition $(p^{[t]}\,i)$. Let  $\pi^{(i-1)}$ be obtained
from $\pi'$ by deleting the last term, which must be $i$. Finally we
set $b_1:=1^{[t]}$, where $\pi^{(1)}(1^{[t]})=1$ for some $t$.

In the above example
$\pi=\bar{5}\bar{\bar{6}}\bar{3}\bar{1}4\bar{\bar{2}}79\bar{\bar{8}}\in\G_{3,9}$
and we successively get
$$\begin{array}{llll}
\pi^{(9)} = \bar{5}\bar{\bar{6}}\bar{3}\bar{1}4\bar{\bar{2}}7{\bf 9}\bar{\bar{8}}, & p=8, & t=0, & b_9=8;\\
\pi^{(8)} = \bar{5}\bar{\bar{6}}\bar{3}\bar{1}4\bar{\bar{2}}7{\bf\bar{\bar{8}}}, & p=8, & t=1, & b_8=\bar{8};\\
\pi^{(7)} = \bar{5}\bar{\bar{6}}\bar{3}\bar{1}4\bar{\bar{2}}{\bf7}, & p=7, & t=0, & b_7=7;\\
\pi^{(6)} = \bar{5}{\bf\bar{\bar{6}}}\bar{3}\bar{1}4\bar{\bar{2}}, & p=2, & t=1, & b_6=\bar{2};\\
\pi^{(5)} = {\bf\bar{5}}\bar{2}\bar{3}\bar{1}4, & p=1, & t=2, & b_5=\bar{\bar{1}};\\
\pi^{(4)} = {\bf\bar{4}}\bar{2}\bar{3}\bar{1}, & p=1, & t=2, & b_4=\bar{\bar{1}};\\
\pi^{(3)} = \bar{\bar{1}}\bar{2}{\bf\bar{3}}, & p=3, & t=2, & b_3=\bar{\bar{3}};\\
\pi^{(2)} = \bar{\bar{1}}{\bf\bar{2}}, & p=2, & t=2, & b_2=\bar{\bar{2}};\\
\pi^{(1)} = {\bf\bar{\bar{1}}}, & p=1, & t=1, & b_1=\bar{1}.\\
\end{array}$$

We can see from the algorithm that the choice of $p^{[t]}$ satisfies
that $i^{[-t]}$ is the $p$-th position in $\pi^{(i)}$. Furthermore,
it also can be seen that $\pi$ can be written as the product
$\prod_{i=1}^{n}(b_i\,i)$, which is a key step for computating
$\sor(\pi)$.


\begin{lem}\label{lem:Bcode_sor}
Suppose that
$b=\text{B-code}(\pi)=(c_1^{[e_1]},c_2^{[e_2]},\ldots,c_n^{[e_n]})$.
Then we have
\begin{equation}\label{eq:sor_formula}
\sor(\pi)=\sum_{i=1}^n \Big(i-c_i+\chi(e_i>0)\cdot\big(2(c_i-1)+e_i\big)\Big)
\end{equation}
and
\begin{equation}\label{eq:ell'_formula}
\ell'(\pi)=n-|\Max^0(b)|.
\end{equation}
\end{lem}

\proof The equality \eqref{eq:sor_formula} can be directly obtained
from \eqref{eq:sor_defi} and we only prove \eqref{eq:ell'_formula}.

Since $\pi=\prod_{i=1}^{n}(b_i\,i)=\prod_{b_i\neq i}(b_i\,i)$, we
have $\ell'(\pi)\leq n-|\Max^0(b)|$. We use induction on $n$ to
prove the equality. If $\pi_n=n$, the assertion is true by taking
$\pi$ as an element in $\G_{r,n-1}$. Assume $\pi_n\neq n$. Let
$\pi(p^{[t]})=n$ for some $p$ and $t$ with $(p,t)\neq (n,0)$. Let
$\pi'=\pi(p^{[t]}\,n)$ and then $\ell'(\pi) = \ell'(\pi')+1$. As
$\pi'$ fixes $n$, we can regard it as an element in $\G_{r,n-1}$ and
by the algorithmic definition of B-code we have
$b':=\text{B-code}(\pi')=(c_1^{[e_1]},\ldots,c_{n-1}^{[e_{n-1}]})$.
Therefore by the induction hypothesis we have
$\ell'(\pi')=(n-1)-|\Max^0(b')|=|\Max^0(b)|$, which completes the
proof. \qed

\begin{lem}\label{lem:Bcode_stirling}
Let $\pi\in\G_{r,n}$ and $b=\text{B-code}(\pi)$. 
Then for each $1\le t \le r-1$ we have
\begin{enumerate}
\item $\Cyc^0(\pi)=\Max^0(b)$ and $\Cyc^t(\pi)=\Max^{r-t}(b)$,
\item $\Lmic^0(\pi)=\Min^0(b)$ and $\Lmic^t(\pi)=\Min^{r-t}(b)$,
\item $\Lmap^0(\pi)=\Rmil^0(b)$ and $\Lmap^t(\pi)=\Rmil^{r-t}(b)$,
\item $\Lmal^0(\pi)=\Rmip^0(b)$ and $\Lmal^t(\pi)=\Rmip^{r-t}(b)$.
\end{enumerate}
\end{lem}
\proof The lemma holds for $n=1$. Assume by induction that the lemma
holds for $n-1$, where $n\geq 2$.

Let $\pi=\pi_1\pi_2\cdots\pi_n$ and $b=(b_1,b_2,\ldots,b_n)$.
By the algorithmic definition of the B-code, there is a colored permutation $\pi'\in\G_{r,n-1}$ such that $b':=\text{B-code}(\pi')=(b_1,b_2,\ldots,b_{n-1})$.
Let $\pi'=C_1C_2\cdots C_k$ be its colored cycle decomposition.
In the following we prove the lemma in two cases according to the position of $n^{[t]}$ in $\pi$.

\emph{{Case 1}.} $\pi_n=n^{[t]}$ for some $0\leq t<r$. In this case
$b=(b_1,b_2,\ldots,b_{n-1},n^{[-t]})$,
$\pi'=\pi_1\pi_2\cdots\pi_{n-1}$ and the colored cycle decomposition
of $\pi$ is $C_1C_2\cdots C_k(n^{[t]})$. It is easy to see that
\begin{enumerate}[(i)]
\item[] $\Cyc(\pi)=\Cyc(\pi')\cup\{n^{[t]}\}$ and $\Max(b)=\Max(b')\cup\{n^{[-t]}\}$,
\item[] $\Lmic(\pi)=\Lmic(\pi')$ and $\Min(b)=\Min(b')$,
\item[] $\Lmap(\pi)=\Lmap(\pi')\cup\{n^{[t]}\}$ and $\Rmil(b)=\Rmil(b')\cup\{n^{[-t]}\}$,
\item[] $\Lmal(\pi)=\Lmal(\pi')\cup\{n^{[t]}\}$ and $\Rmip(b)=\Rmip(b')\cup\{n^{[-t]}\}$,
\end{enumerate}
and the result follows by induction.

\emph{{Case 2}.} $\pi_p=n^{[t]}$ for some $1\leq p<n$ and $0\leq
t<r$. By the definition of wreath product, $\pi(p^{[i]})=n^{[t+i]}$
for any $i$. Then we have $\pi(p^{[-t]})=n$ and thus
$\pi'=\pi\cdot(p^{[-t]}\,n)$. Therefore,
$b=(b_1,b_2,\ldots,b_{n-1},p^{[-t]})$ and hence
\begin{enumerate}[(i)]
\item[] $\Max(b)=\Max(b')$,
\item[] $\Min(b)=\Min(b')\cup\{n^{[-t]}\}$ if $p=1$, and $\Min(b)=\Min(b')$ otherwise,
\item[] $\Rmil(b)=\{i^{[j]}\in\Rmil(b'):\,i<p\}\cup\{p^{[-t]}\}$,
\item[] $\Rmip(b)=\{i^{[j]}\in\Rmip(b'):\,i<n\}\cup\{n^{[-t]}\}$.
\end{enumerate}
Now it suffices to show that
\begin{enumerate}[(i)]
\item $\Cyc(\pi)=\Cyc(\pi')$,
\item $\Lmic(\pi)=\Lmic(\pi')\cup\{n^{[t]}\}$ if $p=1$, and $\Lmic(\pi)=\Lmic(\pi')$ otherwise,
\item $\Lmap(\pi)=\{i^{[j]}\in\Lmap(\pi'):\,i<p\}\cup\{p^{[t]}\}$,
\item $\Lmal(\pi)=\{i^{[j]}\in\Lmal(\pi'):\,i<n\}\cup\{n^{[t]}\}$.
\end{enumerate}

Assume that $\pi_n=m^{[s]}$ for some $1\leq m<n$ and $0\leq s<r$.
From $\pi'=\pi\cdot(p^{[-t]}\,n)$ we obtain
\begin{equation}\label{eq:pi&pi'}
\begin{array}{lccccccccc}
\pi & = & \pi_1 & \cdots & \pi_{p-1} & n^{[t]} & \pi_{p+1} & \cdots & \pi_{n-1} & m^{[s]}\\
\pi' & = & \pi_1 & \cdots & \pi_{p-1} & m^{[s+t]}  & \pi_{p+1} & \cdots & \pi_{n-1} &
\end{array}
\end{equation}

The cases (iii) and (iv) can be obtained directly from
\eqref{eq:pi&pi'}. So we only consider (i) and (ii) in what follows.

(i) Assume that in the colored cycle decomposition of $\pi'$, $C_h$
is the cycle containing the letter $p$, namely, $C_h=(\cdots \,
p^{[t']}\,m^{[s+t]}\,\cdots)$ for some color $t'$. By
\eqref{eq:pi&pi'}, the colored cycle decomposition of $\pi$ must be
$C_1\cdots C_{h-1} \widehat{C_h} C_{h+1}\cdots C_k$ with
$\widehat{C_h}=(\cdots \, p^{[t']}\,n^{[t]}\,m^{[s]}\,\cdots)$.
Notice that the insertion of letter $n$ does not affect the choice
of the smallest letter of $\widehat{C_h}$. Hence
$\Cyc(\pi)=\Cyc(\pi')$, as desired.

(ii) We denote by $g_1\,g_2\,g_3\,\cdots$ the word
$\pi'(1)\pi'^2(1)\pi'^3(1)\cdots$. When $p=1$, it is easy to see
that $g_1=m^{s+t}$ and the word $\pi(1)\pi^2(1)\pi^3(1)\cdots$ will
be $n^{[t]}\,m^{[s+t]}\,g_2\,g_3\,\cdots$. Then we have
$\Lmic(\pi)=\Lmic(\pi')\cup\{n^{[t]}\}$. When $p>1$, the word
$\pi(1)\pi^2(1)\pi^3(1)\cdots$ will be the same as
$g_1\,g_2\,g_3\,\cdots$, or be obtained from it by inserting letters
$n,\bar{n},\ldots$ into some places after $g_1$. Therefore,
$\Lmic(\pi)=\Lmic(\pi')$ in this case. This completes the proof.
\qed

\medskip

Finally, we define $$\phi=(\text{B-code})^{-1}\circ(\text{A-code}).$$

\medskip
\section{Main Theorem A} \label{SectionMainA}

\subsection{Colored permutations} \label{sec:restricted}
We first prove the Main Theorem A in the case of $\G_{r,n}$.

\begin{thm}\label{thm:bijection}
For $\pi\in\G_{r,n}$ we have
\begin{align*}
\big(&\ell,\Rmil^0,\Rmil^1,\ldots,\Rmil^{r-1},\Lmil^0,\Lmil^1,\ldots,\Lmil^{r-1}, \\
&\hspace{1cm} \Lmal^0,\Lmal^1,\ldots,\Lmal^{r-1},\Lmap^0,\Lmap^1,\ldots,\Lmap^{r-1}\big)\pi \\
= \big(&\sor,\Cyc^0,\Cyc^{r-1},\ldots,\Cyc^{1},\Lmic^0,\Lmic^{r-1},\ldots,\Lmic^{1}, \\
&\hspace{1cm} \Lmal^0,\Lmal^{r-1},\ldots,\Lmal^{1},\Lmap^0,\Lmap^{r-1},\ldots,\Lmap^{1}\big)\phi(\pi).
\end{align*}
\end{thm}

\proof The proof is done by combining the Lemma~\ref{lem:Acode_ell},
\ref{lem:Acode_stirling}, \ref{lem:Bcode_sor} and \ref{lem:Bcode_stirling}.
\qed

A direct corollary is the following.
\begin{cor}\label{cor:equidistribution_Grn}
The pair of joint (set-valued) statistics
\begin{align*}
\big(&\ell,\Rmil^0,\Rmil^1,\ldots,\Rmil^{r-1},\Lmil^0,\Lmil^1,\ldots,\Lmil^{r-1}, \\
&\hspace{1cm} \Lmal^0,\Lmal^1,\ldots,\Lmal^{r-1},\Lmap^0,\Lmap^1,\ldots,\Lmap^{r-1}\big)
\end{align*}
and
\begin{align*}
\big(&\sor,\Cyc^0,\Cyc^{r-1},\ldots,\Cyc^{1},\Lmic^0,\Lmic^{r-1},\ldots,\Lmic^{1}, \\
&\hspace{1cm} \Lmal^0,\Lmal^{r-1},\ldots,\Lmal^{1},\Lmap^0,\Lmap^{r-1},\ldots,\Lmap^{1}\big)
\end{align*}
have the same distribution over $\G_{r,n}$.
\end{cor}

It is interesting to see that the superindices of each kind of
statsitics in the first tuple are $0,1,2,\dots, r-1$ while in the
second $0,r-1,r-2,\dots ,1$.

\begin{cor}\label{cor:equidistribution_stirling}
The followings hold.
\begin{enumerate}
\item $\cyc^0,\rmin^0,\lmin^0,\lmax^0$ and $\lmic^0$ have the same distribution with $n-\ell'$ over $\G_{r,n}$.
\item $\cyc,\rmin,\lmin,\lmax$ and $\lmic$ are equidistributed over $\G_{r,n}$.
\end{enumerate}
\end{cor}
\proof
(1) Denote by $\stat_1\sim\stat_2$ the two statistics $\stat_1$ and $\stat_2$ having the same distribution over $\G_{r,n}$.
By \eqref{eq:ell'_formula} and Lemma~\ref{lem:Bcode_stirling}(1), we have $$\ell'(\pi)=n-\cyc^0(\pi),$$
and therefore $(n-\ell')\sim \cyc^0$.

We know from corollary~\ref{cor:equidistribution_Grn} that $\cyc^0\sim\rmin^0$ and $\lmin^0\sim\lmic^0$.
Also from \eqref{eq:pi_inverse_1} we deduce that $\rmin^0\sim\lmax^0$.
Finally, by defining $\pi^r:=\pi(n)\cdots\pi(1)$, the reverse of $\pi$, it is obvious that $\Lmil(\pi)=\Rmil(\pi^r)$
and then we have $\lmin^0\sim\rmin^0$. This completes the proof.

(2) The proof is similar to that of (1) and is omitted.
\qed


\subsection{Colored permutations on a Ferrers shape} \label{sec:restricted}

Let $\ff=(f_1,f_2,\ldots,f_n)$ be a non-decreasing sequence of integers with $1\leq f_1\leq f_2\leq \cdots \leq f_n\leq n$.
The set of \emph{colored restricted permutations} $\G_{r,n,\ff}$ is defined as
$$\G_{r,n,\ff}:=\{\pi=(\sigma,\mathbf{z})\in\G_{r,n}:\,\sigma(i)\leq f_i, 1\leq i\leq n\}.$$
Note that $\G_{r,n,\ff}=\G_{r,n}$ when $\ff=(n,n,\ldots,n)$.
We say that $\ff'$ dominates $\ff$, denoted by $\ff\triangleleft\ff'$, if $f_i\leq f_i'$ for all $i$.
Note that $\G_{r,n,\ff}\subseteq\G_{r,n,\ff'}$ if $\ff\triangleleft\ff'$.

For $\pi\in\G_{r,n}$ define its \emph{minimum sequence} $f(\pi)$ by
$$f(\pi):=\ff \text{ such that } \pi\in\G_{r,n,\ff} \text{ and } \ff\triangleleft\ff' \text{ whenever }\pi\in\G_{r,n,\ff'}.$$
$f(\pi)$ can be easily obtained from $\Lmap(\sigma)$ and $\Lmal(\sigma)$. 
Namely, it is the unique non-decreasing integer sequence $\ff$ such that $\Lmap(\sigma)=\Lmap(\ff)$ and
$\Lmal(\sigma)=\Lmal(\ff)$ by regarding $\ff$ as a word $f_1f_2\cdots f_n$.

More precisely, let $\Lmap(\sigma)=\{i_1,i_2,\ldots,i_s\}$ and $\Lmal(\sigma)=\{\sigma_{i_1},\sigma_{i_2},\ldots,\sigma_{i_s}\}$, where $i_1<\cdots<i_s$.
Then one has $f(\pi)=(f_1,f_2,\ldots,f_n)$ with
$$f_{i_j}=f_{i_j+1}=\cdots=f_{i_{j+1}-1}=\sigma_{i_j}$$
for $j=1,\ldots,s-1$ and $f_{i_s}=\cdots=f_n=\sigma_{i_s}$.

For example, let $\pi=(\sigma,\mathbf{z})=(361475928$, $(3,0,3,2,2,1,0,3,1))\in\G_{4,9}$.
Then $\Lmap(\sigma)=\{1,2,5,7\}$, $\Lmal(\sigma)=\{3,6,7,9\}$ and $f(\pi)=(3,6,6,6,7,7,9,9,9)$.

\begin{lem}\label{lem:restricted}
Let $\pi\in\G_{r,n}$ and $\ff=f(\pi)$.
Then $\phi(\pi)\in\G_{r,n,\ff}$.
\end{lem}
\proof
Let $\pi=(\sigma,\mathbf{z})$ and $\phi(\pi)=(\sigma',\mathbf{z}')$.
By Theorem~\ref{thm:bijection} and the definition of $\Lmap,\Lmal$, we have
\begin{equation}\label{eq:restricted_1}
\Lmap(\sigma)=\bigcup_{t=0}^r\Lmap^t(\pi)=\bigcup_{t=0}^r\Lmap^t(\phi(\pi))=\Lmap(\sigma')
\end{equation}
and
\begin{equation}\label{eq:restricted_2}
\Lmal(\sigma)=\bigcup_{t=0}^r\Lmal^t(\pi)=\bigcup_{t=0}^r\Lmal^t(\phi(\pi))=\Lmal(\sigma').
\end{equation}
Since $\ff$ is determined by $\Lmap$ and $\Lmal$, by \eqref{eq:restricted_1} and \eqref{eq:restricted_2} we have $f(\phi(\sigma))=f(\sigma)=\ff$
and the result follows.
\qed

Finally we come to our main theorem.
\begin{thm}[\textbf{Main Theorem A}] \label{thm:equi_restricted}
Given $r,n,\ff$. Then the pair of joint (set-valued) statistics
\begin{align*}
\big(&\ell,\Rmil^0,\Rmil^1,\ldots,\Rmil^{r-1},\Lmil^0,\Lmil^1,\ldots,\Lmil^{r-1}, \\
&\hspace{1cm} \Lmal^0,\Lmal^1,\ldots,\Lmal^{r-1},\Lmap^0,\Lmap^1,\ldots,\Lmap^{r-1}\big)
\end{align*}
and
\begin{align*}
\big(&\sor,\Cyc^0,\Cyc^{r-1},\ldots,\Cyc^{1},\Lmic^0,\Lmic^{r-1},\ldots,\Lmic^{1}, \\
&\hspace{1cm} \Lmal^0,\Lmal^{r-1},\ldots,\Lmal^{1},\Lmap^0,\Lmap^{r-1},\ldots,\Lmap^{1}\big)
\end{align*}
have the same distribution over $\G_{r,n,\ff}$.
\end{thm}
\proof By Theorem~\ref{thm:bijection}, it suffices to show that if
$\pi\in\G_{r,n,\ff}$, then $\phi(\pi)\in\G_{r,n,\ff}$. And the proof
is done from Lemma~\ref{lem:restricted}.
\qed

\section{Main Theorem B} \label{SectionMainB}
Recall that given $\ff$ we define $H(\ff):=(h_1,h_2,\ldots,h_n)$, where $h_i$ is the smallest possible index at which the letter $i$ can appear in a colored permutation $\sigma\in \G_{r,n,\ff}$. 
It is routine to verify that $$\Rmil(H(\ff))=\Lmap(\ff) \quad \mbox{ and } \quad \Rmip(H(\ff))=\Lmal(\ff)$$ by regarding $\ff$ and $H(\ff)$ as words.

Similar to that of $\G_{r,n}$, all colored restricted permutations of $\G_{r,n,\ff}$ can be generated recursively.
Let $\Psi_1:=1+(\bar{1}\,1)+\cdots+(1^{[r-1]}\,1)$ and for $j\geq 2$
$$\Psi_j:=1+\sum_{i=h_j}^{j-1}(i\,j)+\sum_{t=1}^{r-1}\sum_{i=h_j}^j(i^{[t]}\,j).$$
We omit the proof of the following as it is similar to that of Lemma~\ref{lem:factorization}.

\begin{lem}\label{lem:factorization_restricted}
Given $r,n,\ff$. Then we have
$$\Psi_1\Psi_2\cdots\Psi_n=\sum_{\pi\in\G_{r,n,\ff}}\pi.$$
\end{lem}

Recall the checking function $\xi_{\mathsf{A}}(\omega)$ defined in
the Introduction. Now we arrive at the Main Theorem B.

\begin{thm}[\textbf{Main Theorem B}] \label{thm:GF_restricted}
Given $r,n, \ff$. Let $H(\ff)=(h_1,\ldots,h_n)$.
Then we have
\begin{align*}
&\sum_{\pi\in\G_{r,n,\ff}}q^{\ell(\pi)}\prod_{t=0}^{r-1}\prod_{i\in\Rmil^{-t}(\pi)}x_{t,i}\prod_{i\in\Lmil^{-t}(\pi)}y_{t,i} \\
=& ~~\sum_{\pi\in\G_{r,n,\ff}}q^{\sor(\pi)}\prod_{t=0}^{r-1}\prod_{i\in\Cyc^t(\pi)}x_{t,i}\prod_{i\in\Lmic^t(\pi)}y_{t,i} \\
=& ~~\left.\prod_{j=1}^n \right( x_{0,j}+q+\cdots+q^{j-h_j-1}+ \xi_{h_j=1}(y_{0,j})q^{j-h_j} \\
&\qquad + \left.\sum_{t=1}^{r-1} \Big(x_{r-t,j}q^{2j+t-2}+q^{2j+t-3}+\cdots+q^{j+h_j+t-1}+ \xi_{h_j=1}(y_{r-t,j})q^{j+h_j+t-2} \Big) \right).
\end{align*}
\end{thm}
\proof
We only consider the second equality as the first one is directly from Theorem~\ref{thm:equi_restricted}.

Let $F_n(q,x_{t,i},y_{t,i}:\,0\leq t<r, 1\leq i\leq n)$ denote the desired generating function (which is clearly a polynomial).
Define the linear mapping ${\it\Psi}:\mathbb{Z}(\G_{r,n,\ff})\to\mathbb{Z}(q,x_{t,i},y_{t,i}:\,0\leq t<r, 1\leq i\leq n)$ by
$${\it\Psi}(\pi):=q^{\sor(\pi)}\prod_{t=0}^{r-1}\prod_{i\in\Cyc^t(\pi)}x_{t,i}\prod_{i\in\Lmic^t(\pi)}y_{t,i}.$$
By Lemma~\ref{lem:factorization_restricted} it suffices to show that
$${\it\Psi}(\Psi_1\Psi_2\cdots\Psi_n)=F_n(q,x_{t,i},y_{t,i}:\,0\leq t<r, 1\leq i\leq n).$$ 
We proceed by induction. 
As $\sor(1)=0$, $\sor(1^{[t]})=r-t$ for $0<t<r$, and $\Cyc(1^{[t]})=\Lmic(1^{[t]})=\{1^{[t]}\}$ for $0\leq t<r$, it is easy to see that 
$${\it\Psi}(\Psi_1) = x_{0,1}y_{0,1}+x_{1,1}y_{1,1}q^{r-1}+x_{2,1}t_{2,1}q^{r-2}+\cdots+x_{r-1,1}y_{r-1,1}q.$$
Let $n\geq 2$ and suppose that ${\it\Psi}(\Psi_1\cdots\Psi_{n-1})=F_{n-1}(q,x_{t,i},y_{t,i}:\,0\leq t<r, 1\leq i\leq n-1)$. 
Notice that $\G_{r,n-1,\ff}$ can be identified with the set $\{\pi\in\G_{r,n,\ff}:\,\pi_n=n\}$ in the window notation. 
Given an element $\pi=\pi_1\cdots\pi_{n-1}n$ in this set, we have
$$
\begin{array}{rllll}
\pi\cdot\Psi_n & =\pi_1\cdots\pi_{n-1}n & +~\pi_1\cdots n\,\pi_{n-1} & +~\cdots & +~\pi_1\cdots\pi_{h_n-1}\,n\,\pi_{h_n+1}\cdots\pi_{n-1}\pi_{h_n} \\ \\
 & +~\pi_1\cdots\pi_{n-1}n^{[r-1]} & +~\pi_1\cdots n^{[r-1]}\pi_{n-1}^{[1]} & +~\cdots & +~\pi_1\cdots\pi_{h_n-1}n^{[r-1]}\pi_{h_n+1}\cdots\pi_{n-1}\pi_{h_n}^{[1]} \\ \\
 & +~\cdots\cdots\cdots & & & \\ \\
 & +~\pi_1\cdots\pi_{n-1}n^{[1]} & +~\pi_1\cdots n^{[1]}\pi_{n-1}^{[r-1]} & +~\cdots & +~\pi_1\cdots\pi_{h_n-1}n^{[1]}\pi_{h_n+1}\cdots\pi_{n-1}\pi_{h_n}^{[r-1]}.
\end{array}
$$
Denote by $\pi'$ any one of the summands above.
Without loss of generality, let the letter $n^{[t]}$, for some $t$, be at the $i$-th position in $\pi'$.
That is, $\pi'=\pi\cdot(i^{[-t]}\,n)$.
Then $\sor(\pi')=\sor(\pi)+n-i$ if $t=0$ and $\sor(\pi')=\sor(\pi)+n+i+r-t-2$ if $t>1$.
Moreover, from the proof of Lemma~\ref{lem:Bcode_stirling}, we have
\begin{align*}
\Cyc(\pi')=
\begin{cases}
\Cyc(\pi) \cup \{n^{[t]}\} & \text{if } \pi'=\pi_1\cdots\pi_{n-1}\,n^{[t]} \text{ for some }t, \\
\Cyc(\pi) & \text{otherwise; }
\end{cases}
\end{align*}
and
\begin{align*}
\Lmic(\pi')=
\begin{cases}
\Lmic(\pi) \cup \{n^{[t]}\} & \text{if } \pi'=n^{[t]}\pi_2\cdots\pi_{n-1}\pi_1^{[-t]} \text{ for some }t, \\
\Lmic(\pi) & \text{otherwise. }
\end{cases}
\end{align*}
Therefore,
\begin{align*}
{\it\Psi}(\pi\cdot\Psi_{n}) =& {\it\Psi}(\pi) \Big( x_{0,n}+q+\cdots+q^{n-h_n-1}+ \xi_{h_n=1}(y_{0,n})q^{n-h_n} \\
\qquad &+ \sum_{t=1}^{r-1} \Big(x_{r-t,n}q^{2n+t-2}+q^{2n+t-3}+\cdots+q^{n+h_n+t-1}+ \xi_{h_n=1}(y_{r-t,n})q^{n+h_n+t-2} \Big) \Big).
\end{align*}

Now it suffices to show that ${\it\Psi}(\Psi_1\cdots\Psi_{n-1}){\it\Psi}(\Psi_n) = {\it\Psi}(\Psi_1\cdots\Psi_{n-1}\Psi_n)$.
This can be done by the same argument in the proof of Theorem~\ref{thm:sor_GF}.
Hence we are done.
\qed

\medskip

We obtain the following corollary by replacing $x_{t,i}$ with $x_t$ and $y_{t,i}$ with $y_t$, for each $0\leq t<r$ and $1\leq i\leq n$.

\begin{cor} \label{cor:GF_restricted}
Given $r,n,\ff$. Let $H(\ff)=(h_1,\ldots,h_n)$. Then we have
\begin{align*}
&\sum_{\pi\in\G_{r,n,\ff}}q^{\ell(\pi)}\prod_{t=0}^{r-1} x_t^{\rmin^{-t}(\pi)} y_t^{\lmin^{-t}(\pi)}
=\sum_{\pi\in\G_{r,n,\ff}}q^{\sor(\pi)}\prod_{t=0}^{r-1} x_t^{\cyc^{t}(\pi)} y_t^{\lmic^{t}(\pi)}\\
=& ~~\left.\prod_{j=1}^n \right( x_{0}+q+\cdots+q^{j-h_j-1}+ \xi_{h_j=1}(y_0)q^{j-h_j} \\
&\qquad + \left.\sum_{t=1}^{r-1} \Big(x_{r-t}q^{2j+t-2}+q^{2j+t-3}+\cdots+q^{j+h_j+t-1}+ \xi_{h_j=1}(y_{r-t})q^{j+h_j+t-2}\Big) \right).
\end{align*}
In particular,
\begin{align*}
&\sum_{\pi\in\G_{r,n}}q^{\ell(\pi)}\prod_{t=0}^{r-1} x_t^{\rmin^{-t}(\pi)} y_t^{\lmin^{-t}(\pi)}
=\sum_{\pi\in\G_{r,n}}q^{\sor(\pi)}\prod_{t=0}^{r-1} x_t^{\cyc^{t}(\pi)} y_t^{\lmic^{t}(\pi)}\\
=& ~~\prod_{j=1}^n \left( x_0+y_0q^{j-1} +\sum_{t=1}^{r-1}q^{j+r-t-1}\left(x_tq^{j-1}+y_t\right) +q[j-2]_q\left(1+q^j[r-1]_q\right) \right).
\end{align*}
\end{cor}

\medskip

Since $\ell'=n-\cyc^0$, we have the following.

\begin{cor}\label{cor:distribution_ell'}
We have
$$\sum_{\pi\in\G_{r,n}}t^{\cyc^0(\pi)} = \prod_{i=1}^n(t+ri-1) \quad \text{and} \quad \sum_{\pi\in\G_{r,n}}t^{\ell'(\pi)} = \prod_{i=1}^n(1+(ri-1)t).$$
\end{cor}


\section{Even-signed permutation group} \label{sec:Dn}

We turn to the case of even-signed permutation group, defined as the subgroup $\D_n$ of $\G_{2,n}$ consisting of those signed permutations $\pi$ with even number of negatives in the window notation $\pi=\pi_1,\cdots,\pi_n$. 
Here we adopt the convention $\bar{i}=-i$.

\subsection{Sorting index}
It is known that $\D_n$ is a Coxeter group generated by $$\mathcal{S}_n^D=\{s^D_0,s_1,\ldots,s_{n-1}\},$$ where $s^D_0$ is the transposition $(\bar{1}\,2)$ and $s_i$ is the transposition $(i\,i+1)$ for $i\geq 1$. 
Let $\ell_D$ be the length function of $\D_n$ with respect to $\mathcal{S}_n^D$.

$\D_n$ can also be generated by
$$\mathcal{T}_n^D:=\{t^D_{ij}:\,1\leq |i|<j\leq
n\}\cup\{t^D_{\bar{i}i}:\,1<i\leq n\},$$ where $t^D_{ij}=(i\,j)$ for
$1\leq |i|<j\leq n$ and $t^D_{\bar{i}i}=(\bar{1}\,1)(\bar{i}\,i)$
for $1<i\leq n$. For $\pi\in\D_n$, let $\tilde{\ell}_D'(\pi)$ be the
minimum number of elements in $\mathcal{T}_n^D$ needed to express
$\pi$. Again we note that $\tilde{\ell}_D'(\pi)$ is not the
reflection length of $\pi$~\cite{Petersen_11}.

Any $\pi\in\D_n$ has a unique factorization in the form
$$\pi=t^D_{i_1 j_1}t^D_{i_2 j_2}\cdots t^D_{i_k j_k}$$ with
$1<j_1<j_2<\cdots<j_k\leq n$. Petersen~\cite{Petersen_11} defined
the sorting index
$$\sor_D(\pi)=\sum_{r=1}^k (j_r-i_r-2\chi(i_r<0)).$$
and proved that it is Mahonian. For example,
$\pi=\bar{3}24\bar{5}1=t^D_{\bar{1}3}t^D_{34}t^D_{\bar{4}5}$ has the
sorting index $\sor_D(\pi)= (3-(-1)-2)+ (4-3)+ (5-(-4)-2)= 10$.

\subsection{Set-valued Stirling statistics}
For $\pi\in\D_n$ we define the (set-valued) statistics
$$\Cyc_D(\pi), \Rmil_D(\pi), \Rmip_D(\pi), \Lmal_D(\pi),
\Lmap_D(\pi), \Lmil_D(\pi), \Lmic_D(\pi)$$ and $$\cyc_D(\pi),
\rmin_D(\pi), \lmax_D(\pi), \lmin_D(\pi), \lmic_D(\pi)$$ by viewing
$\pi$ as an element in $\G_{2,n}$. The subscript $D$ is to emphasize
that these statistics are considered in $\D_n$. We also define some
new `twisted' statistics:
\begin{enumerate}
\item $\Cyc_D^+$, the set of \emph{twisted balanced cycles}: $$\Cyc_D^{+}(\pi):=\Cyc_D^{0}(\pi)\cup\{1\}.$$
\item $\Cyc_D^-$, the set of \emph{twisted unbalanced cycles}: $$\Cyc_D^{-}(\pi):=\Cyc_D^{1}(\pi)\setminus\{1\}.$$
\item $\Rmil_D^{+}$, the set of \emph{twisted positive right-to-left minimum letters}: $$\Rmil_D^{+}(\pi):=\Rmil_D^{0}(\pi)\cup\{1\}.$$
\item $\Rmil_D^{-}$, the set of \emph{twisted negative right-to-left minimum letters}: $$\Rmil_D^{-}(\pi):= \Rmil_D^{1}(\pi)\setminus\{1\}.$$
\end{enumerate}
Denote by $\cyc_D^+,\cyc_D^-,\rmin_D^+,\rmin_D^-$ their cardinalities respectively.

\subsection{A bijection}
The key ingredient in this section is the bijection
$\psi:\D_n\to\D_n$, introduced by Chen-Guo-Gone~\cite{Chen_13},
which is the composition
$$\psi=(\text{D-code})^{-1}\circ (\text{C-code})$$
 of the C-code and the D-code on $\D_n$.

\medskip
\noindent \textbf{Algorithm for C-code}. 
For $\pi\in\D_n$ we construct a sequence of $n$ even-signed permutations $\pi=\pi^{(n)},\pi^{(n-1)},\ldots,\pi^{(1)}$ such that $\pi^{(i)}\in\D_i$ and meanwhile obtain the C-code $(c_1,c_2,\dots,c_n)$. 
For $i$ from $n$ down to $2$ we consider the letter $i$ or $\bar{i}$ in $\pi^{(i)}$. If $i$ appears at the $p$-th position in $\pi^{(i)}$, then we define $c_i=p$ and let $\pi^{(i-1)}$ be obtained from $\pi^{(i)}$ by deleting the letter $i$. 
If $\bar{i}$ appears at the $p$-th position in $\pi^{(i)}$, then we first define $c_i=-p$, let $\pi'$ be obtained by deleting $\bar{i}$ and then obtain $\pi^{(i-1)}$ from $\pi'$ by changing the sign of the first letter. 
Finally, set $c_1:=\pi^{(1)}(1)$. 
It is easy to see that $\pi^{(1)}$ is always the identity permutation and hence $c_1=1$.

For example, let $\pi=\bar{5}\bar{2}\bar{1}\bar{3}4$, then
$$\begin{array}{ll}
\pi^{(5)} = {\bf\bar{5}}\,\bar{2}\,\bar{1}\,\bar{3}\,4, & c_5=-1, \\
\pi^{(4)} = 2\,\bar{1}\,\bar{3}\,{\bf 4}, & c_4=4, \\
\pi^{(3)} = 2\,\bar{1}\,{\bf\bar{3}}, & c_3=-3, \\
\pi^{(2)} = {\bf\bar{2}}\,\bar{1}, & c_2=-1, \\
\pi^{(1)} = {\bf 1}, & c_1=1, \\
\end{array}$$
hence $\text{C-code}(\pi)=(1,-1,-3,4,-1)$.

Let $\mathsf{SE}_n^D$ be the set of integer sequences $(c_1,c_2,\ldots,c_n)$ such that $c_1=1$ and for $i\geq 2$, $c_i\in [-i,i]\backslash\{0\}$.
It is obvious that C-code is a bijection from $\D_n$ onto $\mathsf{SE}_n^D$.

\begin{lem}\label{lem:Dn_Ccode}
For $\pi\in\D_n$ we have
\begin{enumerate}
\item $(\Lmil_D,\Lmap_D,\Lmal_D)\pi=(\Min,\Rmil_D,\Rmip_D)\text{C-code}(\pi)$.
\item $(\Rmil_D^+,\Rmil_D^-)\pi = (\Max^0,\Max^1)\text{C-code}(\pi)$.
\end{enumerate}
\end{lem}
\proof
(1) Let $a=\text{A-code}(\pi)$ and $c=\text{C-code}(\pi)$.
By definition one has $$|a_i|=|c_i| \quad \text{for } i=1,2,\ldots,n.$$
Hence the result follows from Lemma~\ref{lem:Acode_stirling}(2) -- (4).

(2) Observe that the sign of each letter on the right of $1$ or $\bar{1}$ will not change during the construction of the C-code.
Assume that $\pi_p=1$ or $\bar{1}$ for some $p$.
Then we have $a_j=c_j$ for $j\in\{\pi_i:\,i>p\}$.
Thus, $a_j=c_j$ for $j\in\Rmil_D(\pi)$.
Since $c_1$ is always $1$, by Lemma~\ref{lem:Acode_stirling}(1), we have
$$\Rmil_D^{+}(\pi)=\Rmil_D^0(\pi)\cup\{1\}={\Max}^0(a)\cup\{1\}={\Max}^0(c)$$ and
$$\Rmil_D^{-}(\pi)=\Rmil_D^0(\pi)\setminus\{1\}={\Max}^0(a)\setminus\{1\}={\Max}^1(c),$$ hence the lemma is proved.
\qed

\medskip

The D-code is defined as follows.

\noindent \textbf{Algorithm for D-code}. For $\pi\in\D_n$, we
generate a sequence of even-signed permutations
$\pi=\pi^{(n)},\pi^{(n-1)},\ldots,\pi^{(1)}$ such that
$\pi^{(i)}\in\D_i$ and meanwhile construct the D-code
$(d_1,d_2,\ldots,d_n)$. For $i$ from $n$ down to $2$ we consider the
letter $i$ or $\bar{i}$ in $\pi^{(i)}$. If $i$ appears at the $p$-th
position, we set $d_i=p$ and let $\pi^{(i-1)}$ be the first $i-1$
terms of $\pi^{(i)}t^D_{p\,i}$. If $\bar{i}$ appears at the $p$-th
position, then we set $d_i=-i$ and let $\sigma^{(i-1)}$ be the first
$i-1$ terms of $\pi^{(i)}t^D_{\bar{p}\,i}$. It can be seen that
$\pi^{(1)}$ is always the identity $1$ and hence $d_1=1$.

For example, if $\pi=\bar{5}\bar{1}\bar{3}4\bar{2}$, then we have
$$\begin{array}{ll}
\pi^{(5)} = {\bf\bar{5}}\,\bar{1}\,\bar{3}\,4\,\bar{2}, & d_5=-1 \\
\pi^{(4)} = 2\,\bar{1}\,\bar{3}\,{\bf 4}, & d_4=4 \\
\pi^{(3)} = 2\,\bar{1}\,{\bf\bar{3}}, & d_3=-3 \\
\pi^{(2)} = {\bf\bar{2}}\,\bar{1}, & d_2=-1 \\
\pi^{(1)} = {\bf 1}, & d_1=1
\end{array}$$
and thus $\text{D-code}(\pi)=(1,-1,-3,4,-1)$.

\begin{lem}\label{lem:Dn_Dcode}
For $\pi\in\D_n$ we have
\begin{enumerate}
\item $(\Lmic_D,\Lmap_D,\Lmal_D)\pi = (\Min,\Rmil_D,\Rmip_D)\text{D-code}(\pi)$.
\item $(\Cyc_D^+,\Cyc_D^-)\pi = (\Max^0,\Max^1)\text{D-code}(\pi)$.
\end{enumerate}
\end{lem}
\proof
(1) Let $b=\text{B-code}(\pi)$ and $d=\text{D-code}(\pi)$.
By definition one has $$|b_i|=|d_i| \quad \text{for } i=1,2,\ldots,n.$$
The result follows from Lemma~\ref{lem:Bcode_stirling}(2) -- (4).

(2) Fix an integer $1<i\leq n$. 
Let $\sigma$ and $\rho$ denote the permutations $\pi^{(i)}$ during the construction of the B-code and
D-code respectively. 
One can see that $|\sigma_1|=|\rho_1|$ and $\sigma_j=\rho_j$ for $j\geq 2$. 
That is, $d_j$ must be equal to $b_j$ whenever $|d_j|\neq 1$. 
Since $d_1$ is always $1$, by Lemma~\ref{lem:Bcode_stirling}(1) we have
$$\Cyc_D^{+}(\pi)=\Cyc_D^0(\pi)\cup\{1\}={\Max}^0(b)\cup\{1\}={\Max}^0(d)$$ and
$$\Cyc_D^{-}(\pi)=\Cyc_D^1(\pi)\setminus\{1\}={\Max}^1(b)\setminus\{1\}={\Max}^1(d),$$ and we are done.
\qed

Note that in \cite[Proposition 4.5]{Chen_13} it is showed that $\tilde{\ell}_D'(\pi)=n-\sum_{i=1}^n \chi(d_i=i)$. 
Hence from Lemma~\ref{lem:Dn_Dcode} we have
\begin{equation}
\tilde{\ell}_D'(\pi)=n-\cyc_D^{+}(\pi).
\end{equation}

\medskip
The next result extends the type $D$ main result
in~\cite{Poznanovic_14}.
\begin{thm}\label{thm:equi_Dn}
For $\pi\in\D_n$ we have $$(\ell_D,\Rmil_D^{+},\Rmil_D^{-},\Lmil_D,\Lmap_D)\pi = (\sor_D,\Cyc_D^+,\Cyc_D^-,\Lmic_D, \Lmap_D)\psi(\pi).$$
\end{thm}
\proof It was showed in \cite{Chen_13} that
$\inv_D(\pi)=\sor_D(\psi(\pi))$ for $\pi\in\D_n$. The theorem is
proved by combining it with Lemma~\ref{lem:Dn_Ccode} and
Lemma~\ref{lem:Dn_Dcode}. \qed

\subsection{Even-signed permutations on a Ferrers shape} We can also extend the result on the permutations restricted to a Ferrers shape.
Similar to the case of $\G_{r,n,\ff}$, for a given integer sequence
$\ff=(f_1,f_2,\ldots,f_n)$ with $1\leq f_1\leq f_2\leq \cdots\leq
f_n\leq n$ we define the set of restricted even-signed permutations
by
$$\D_{n,\ff} := \{\pi\in\D_n:\,|\pi_i|\leq f_i, 1\leq i\leq n\}.$$
The minimum sequence $f(\pi)$ is similarly determined by both
$\Lmap(|\pi|)$ and $\Lmal(|\pi|)$, where $|\pi|$ is the permutation
$|\pi_1|\cdots|\pi_n|\in\mathfrak{S}_n$. We have the following.

\begin{lem}\label{lem:restricted_Dn}
Let $\pi\in\D_n$ and $f(\pi)=\ff$. Then $\psi(\pi)\in\D_{n,\ff}$.
\end{lem}
\proof
From the proof of Lemma~\ref{lem:restricted} we have $$(\Lmap,\Lmal)|\pi| = (\Lmap,\Lmal)\psi(|\pi|).$$
Therefore, $f(\psi(\pi))=\ff$ and thus the result follows.
\qed

By combining Theorem~\ref{thm:equi_Dn} and
Lemma~\ref{lem:restricted_Dn}, we obtain the first main result of
this section.

\begin{thm}\label{thm:equi_restricted_Dn}
Given $n$ and $\ff$. Then the pair of (set-valued) statistics
$$(\ell_D,\Rmil_D^{+},\Rmil_D^{-},\Lmil_D,\Lmap_D,\Lmal_D) \quad\mbox{and}\quad
(\sor_D,\Cyc_D^+,\Cyc_D^-,\Lmic_D,\Lmap_D,\Lmal_D)$$ have the same
joint distribution over $\D_{n,\ff}$.
\end{thm}

Now we look at the generating function. 
Define elements $\Theta_1,\Theta_2,\ldots,\Theta_n$ of the group algebra of $\D_{n,\ff}$ by $\Theta_1:=1$ and for $j\geq 2$ 
$$\Theta_j:=1 + \sum_{i=h_j}^{j-1}t^D_{ij} + \sum_{i=h_j}^{j}t^D_{\bar{i}j},$$ 
where $H(\ff)=(h_1,\ldots,h_n)$.
The next lemma is derived by a similar argument in the proof of Lemma~\ref{lem:factorization}.

\begin{lem}\label{lem:factorization_Dn}
We have
$$\Theta_1\Theta_2\cdots\Theta_n=\sum_{\pi\in\D_{n,\ff}}\pi.$$
\end{lem}

Our second main result is a generating function over $\D_{n,\ff}$, analogous to Theorem~\ref{thm:GF_restricted}.

\begin{thm}\label{thm:GF_Dn}
Given $n$ and $\ff$ with $H(\ff)=(h_1,\ldots,h_n)$, we have
\begin{align*}
&\sum_{\pi\in\D_{n,\ff}}q^{\ell_D(\pi)}u^{\lmin_D(\pi)} \prod_{i\in\Rmil_D^+(\pi)}t_i \prod_{i\in\Rmil_D^-(\pi)}s_i \\
=& ~~\sum_{\pi\in\D_{n,\ff}}q^{\sor_D(\pi)}u^{\lmic_D(\pi)} \prod_{i\in\Cyc_D^+(\pi)}t_i \prod_{i\in\Cyc_D^-(\pi)}s_i \\
=& ~~t_1u \prod_{j=2}^n \left( t_j+q+\cdots+q^{j-h_j-1}+ \xi_{h_j=1}(2u)q^{j-h_j} +q^{j+h_j-1}+\cdots+q^{2j-3}+s_jq^{2j-2}\right).
\end{align*}
\end{thm}
\proof The first equality follows from
Theorem~\ref{thm:equi_restricted_Dn}. For the second one, let
$F_n(q,u,t_i,s_i:\,1\leq i\leq n)$ denote the desired generating
function. Define the linear mapping
$\theta:\mathbb{Z}(\D_{n,\ff})\to\mathbb{Z}(q,u,t_i,s_i:\,1\leq
i\leq n)$ by
$$\theta(\pi):=q^{\sor_D(\pi)}u^{\lmic_D(\pi)} \prod_{i\in\Cyc_D^+(\pi)}t_i \prod_{i\in\Cyc_D^-(\pi)}s_i.$$
By Lemma~\ref{lem:factorization_Dn} it suffices to show that $$\theta(\Theta_1\Theta_2\cdots\Theta_n)=F_n(q,u,t_i,s_i:\,1\leq i\leq n).$$
We proceed by induction.
It is easy to see that $\theta(\Theta_1)=t_1u$ and
\begin{equation*}
\theta(\Theta_1\Theta_2)=
\begin{cases}
t_1u(t_2+s_2q^2) \mbox{ if } \ff=(1,2) \mbox{ or } \\
t_1u(t_2+2uq+s_2q^2) \mbox{ if } \ff=(2,2).
\end{cases}
\end{equation*}

Let $n\geq 3$ and suppose that
$\theta(\Theta_1\cdots\Theta_{n-1})=F_{n-1}(q,u,t_i,s_i:\,1\leq
i\leq n-1)$. Similar to the proof of
Theorem~\ref{thm:GF_restricted}, we identify elements of
$\D_{n-1,\ff}$ with the set $\{\pi\in\D_{n,\ff}:\,\pi_n=n\}$. Given
an element $\pi=\pi_1\cdots\pi_{n-1}n$ in this set, we have
$$
\begin{array}{rllll}
\pi\cdot\Theta_n & =\pi_1\pi_2\cdots\pi_{n-1}n & +~\pi_1\cdots n\,\pi_{n-1} & +~\cdots & +~\pi_1\cdots\pi_{h_n-1}\,n\,\pi_{h_n+1}\cdots\pi_{n-1}\pi_{h_n} \\
 & +~\overline{\pi_1}\pi_2\cdots\pi_{n-1}\bar{n} & +~\pi_1\cdots \bar{n}\overline{\pi_{n-1}} & +~\cdots & +~\pi_1\cdots\pi_{h_n-1}\bar{n}\pi_{h_n+1}\cdots\pi_{n-1}\overline{\pi_{h_n}}.
\end{array}
$$
Denote by $\pi'$ any one of the summands above.
Without loss of generality, let the letter $n$ or $\bar{n}$ be at the $i$-th position in $\pi'$.
That is, $\pi'=\pi t^D_{in}$ or $\pi t^D_{\bar{i}n}$.
Then $\sor_D(\pi')=\sor_D(\pi)+n-i$ if $n$ has a positive sign and $\sor_D(\pi')=\sor_D(\pi)+n+i-2$ otherwise.
Since $1$ is always counted in $\Cyc_D^+$, by the proof of Lemma~\ref{lem:Bcode_stirling}, we have
\begin{align*}
\Cyc_D^+(\pi')=\begin{cases}
\Cyc_D^+(\pi) \cup \{n\} & \text{if } \pi'=\pi_1\cdots\pi_{n-1}\,n, \\
\Cyc_D^+(\pi) & \text{otherwise.}
\end{cases}
\end{align*}
and
\begin{align*}
\Cyc_D^-(\pi')=\begin{cases}
\Cyc_D^-(\pi) \cup \{n\} & \text{if } \pi'=\pi_1\cdots\pi_{n-1}\,\bar{n}, \\
\Cyc_D^-(\pi) & \text{otherwise.}
\end{cases}
\end{align*}
Moreover, it is not hard to verify that
\begin{align*}
\lmic_D(\pi')=\begin{cases}
\lmic_D(\pi)+1 & \text{if } \pi'=n\pi_2\cdots\pi_{n-1}\pi_1 \text{ or } \bar{n}\pi_2\cdots\pi_{n-1}\overline{\pi_1}\\
\lmic_D(\pi) & \text{otherwise.}
\end{cases}
\end{align*}
Thus, we have
\begin{align*}
\theta &(\pi\cdot\Theta_n) = \\
&\theta(\pi) \left(t_n+q+\cdots+q^{n-h_n-1}+ \xi_{h_n=1}(2u)q^{n-h_n} +q^{n+h_n-1}+\cdots+q^{2n-3}+s_nq^{2n-2}\right)
\end{align*}
and therefore
\begin{align*}
\theta &(\Theta_1\cdots\Theta_{n-1}\Theta_{n}) = \theta\left(\sum_{\pi\in \D_{n,\ff},\,\pi_n=n} \pi\cdot\Theta_n \right) = \sum_{\pi\in \D_{n,\ff},\,\pi_n=n} \theta(\pi\cdot\Theta_n)\\
&= \left(t_n+q+\cdots+q^{n-h_n-1}+ \xi_{h_n=1}(2u)q^{n-h_n} +q^{n+h_n-1}+\cdots+q^{2n-3}+s_nq^{2n-2}\right) \sum_{\pi\in\D_{n-1,\ff}}\theta(\pi)\\
&= \left(t_n+q+\cdots+q^{n-h_n-1}+ \xi_{h_n=1}(2u)q^{n-h_n} +q^{n+h_n-1}+\cdots+q^{2n-3}+s_nq^{2n-2}\right) \\
& \hspace{9cm} \cdot F_{n-1}(q,u,t_i,s_i:\,1\leq i\leq n-1)\\
&= F_n(q,u,t_i,s_i:\,1\leq i\leq n).
\end{align*}
\qed

 By replacing $t_i$ with $t$ and $s_i$ with $s$ for all $i$,
 we obtain the following result.

\begin{cor} \label{cor:GF_Dn}
Given $n$ and $\ff$ with $H(\ff)=(h_1,\ldots,h_n)$, we have
\begin{align*}
&\sum_{\pi\in\D_{n,\ff}}q^{\ell_D(\pi)} u^{\lmin_D(\pi)} t^{\rmin_D^+(\pi)} s^{\rmin_D^-(\pi)}
=\sum_{\pi\in\D_{n,\ff}}q^{\sor_D(\pi)} u^{\lmic_D(\pi)} t^{\cyc_D^+(\pi)} s^{\cyc_D^-(\pi)} \\
=& ~~tu \prod_{j=2}^n \left(t+q+\cdots+q^{j-h_j-1}+ \xi_{h_j=1}(2u)q^{j-h_j} +q^{j+h_j-1}+\cdots+q^{2j-3}+sq^{2j-2}\right).
\end{align*}
In particular,
\begin{align*}
&\sum_{\pi\in\D_{n}}q^{\ell_D(\pi)} u^{\lmin_D(\pi)} t^{\rmin_D^+(\pi)} s^{\rmin_D^-(\pi)}
=\sum_{\pi\in\D_{n}}q^{\sor_D(\pi)} u^{\lmic_D(\pi)} t^{\cyc_D^+(\pi)} s^{\cyc_D^-(\pi)} \\
=& ~~tu \prod_{j=2}^n \left(t+q+\cdots+q^{j-2}+2uq^{j-1}+q^{j}+\cdots+q^{2j-3}+sq^{2j-2}\right).
\end{align*}
\end{cor}

\medskip

Since $\tilde{\ell}_D'=n-\cyc_D^{+}$, we have the following.

\begin{cor}\label{cor:distribution_ell'_Dn}
We have
$$\sum_{\pi\in\D_{n}}t^{\cyc^+(\pi)} = t\prod_{i=2}^n(t+2i-1) \quad \text{and} \quad \sum_{\pi\in\D_{n}}t^{\tilde{\ell}_D'(\pi)} = \prod_{i=2}^n(1+(2i-1)t).$$
\end{cor}

\section{Concluding Remark}\label{sec:conclusion}
 In the notation of complex reflection groups, $\G_{r,n}$ is denoted by $G(r,1,n)$ and $\D_n$ by $G(2,2,n)$. Hence it is natural to ask if one can have the
sorting index and analogous (set-valued) equidistribution results on
$G(r,2,n)$, or even better, $G(r,r,n)$ or $G(r,p,n)$. We leave these
questions to the interested readers.

\rm


\end{document}